\documentclass[12pt,final,leqno]{amsart}
\RequirePackage{lmodern} %Not1PartOfTexFile

\usepackage{amsfonts}
\usepackage{amsmath}
\usepackage{amssymb}
\usepackage{amsthm}
\usepackage{appendix}
\usepackage{changepage}
\usepackage{enumerate}
\usepackage{enumitem}
\usepackage{eucal}
\usepackage{latexsym}
\usepackage{mathrsfs}
\usepackage{scrextend}
\usepackage{textcomp}
\usepackage{todonotes}
\usepackage{verbatim}
\usepackage{pdflscape}
\usepackage{lscape}

\relax							

\newlength{\pardefault}							
\newenvironment{psmallmatrix}
  {\left(\begin{smallmatrix}}
  {\end{smallmatrix}\right)}

\newcommand{\innerProductReg}[2]{\left\langle #1 , #2 \right\rangle}

\newcommand{\innerProductIndefiniteReg}[2]{\left [ #1 , #2 \right]}
\newcommand{\innerProductIndefiniteTri}[3]{\left [ #1 , #2 \right]_{#3}}
\newcommand{\innerProductIndefiniteQuad}[4]{\left [ #1 , #2 \right]_{#3}^{#4}}

\newcommand{\norm}[1]{\| #1  \|}

\newcommand{\normA}[1]{ \left| \left| { #1 } \right| \right| }
\newcommand{\normTreA}[3]{ \left| \left| { #1 } \right| \right|  _{#2}^{#3} }

\newcommand{\bracketsA}[1]{\left( #1 \right)}

\newcommand{\abs}[1]{\left| #1 \right|}

\newcommand{\real}{\mathfrak{Re} ~}
\newcommand{\imag}{\mathfrak{Im} ~}

\newcommand{\restrict}[1]{\raisebox{-.5ex}{$|$}_{#1}}
\newcommand\restr[2]{{% we make the whole thing an ordinary symbol
  \left.\kern-\nulldelimiterspace % automatically resize the bar with \right
  #1 % the function
  \vphantom{\big|} % pretend it's a little taller at normal size
  \right|_{#2} % this is the delimiter
  }}
\newcommand{\defEq}{\overset{\text{def} } {=}}

\newcommand{\forceOddPage}{
\ifodd\theCurrentPage\pagebreak\fi}
\newcommand{\forceEvenPage}{
\ifodd\value{\theCurrentPage}\newpage\else\pagebreak\fi}
\newcommand{\ignore}[1]{}

\newtheorem{Pa}{Paper}[section]
\newtheorem{Tm}[Pa]{{\bf Theorem}}

\newtheorem{Cy}[Pa]{{\bf Corollary}}
\newtheorem{Rk}[Pa]{{\bf Remark}}

\newtheorem{theorem}[Pa]{{\bf Theorem}}
\newtheorem{lem}[Pa]{{\bf Lemma}}

\newtheorem{definition}[Pa]{{\bf Definition}}

\newtheorem{Step}[]{{\bf Step}}

\newenvironment{pf}[1][\unskip]{
\par
\noindent
\paragraph{{\bf Proof #1:}}
\noindent
}
{\hfill$\square$\\}

\usepackage{hyperref}
\hypersetup{
  colorlinks,
  citecolor	= violet,
  linkcolor	= blue,
  urlcolor	= red}

% ...

\begin{document}

%%%%%%%%%%%%%%%%%%%%%%%%%%%%%%%%%%%%%
%									%
%			Title					%
%									%
%%%%%%%%%%%%%%%%%%%%%%%%%%%%%%%%%%%%%
\title
[Commuting operators over Pontryagin spaces]
{Commuting operators over Pontryagin spaces with applications to system theory}

%%%%%%%%%%%%%%%%%%%%%%%%%%%%%%%%%%%%%
%									%
%			Authors					%
%									%
%%%%%%%%%%%%%%%%%%%%%%%%%%%%%%%%%%%%%
\author[D. Alpay]{Daniel Alpay}
\address{(DA) Schmid College of Science and Technology,
Chapman University, One University Drive Orange, California 92866, USA}
\email{alpay@chapman.edu}

\author[A. Pinhas]{Ariel Pinhas}
\address{(AP) Department of mathematics,
Ben-Gurion University of the Negev, P.O. Box
653, Beer-Sheva 84105, Israel}
\email{arielp@post.bgu.ac.il}

\author[V. Vinnikov]{Victor Vinnikov}
\address{(VV) Department of mathematics,
Ben-Gurion University of the Negev, P.O. Box
653, Beer-Sheva 84105, Israel}
\email{vinnikov@math.bgu.ac.il}

\date{}

\thanks{The first author thanks the Foster G. and Mary McGaw Professorship in
Mathematical Sciences, which supported this research}

%%%%%%%%%%%%%%%%%%%%%%%%%%%%%%%%%%%%%
%									%
%			Abstract				%
%									%
%%%%%%%%%%%%%%%%%%%%%%%%%%%%%%%%%%%%%
\begin{abstract}

In this paper we extend vessel theory, or equivalently, the theory of overdetermined $2D$ systems to the Pontryagin space setting.
We focus on realization theorems of the various characteristic functions associated to such vessels.
In particular, we develop an indefinite version of de Branges-Rovnyak theory over real compact Riemann surfaces.
To do so, we use the theory of contractions in Pontryagin spaces and the theory of analytic kernels with a finite number of negative squares.
Finally, we utilize the indefinite de Branges-Rovnyak theory on compact Riemann surfaces
in order to prove a Beurling type theorem on indefinite Hardy spaces on finite bordered Riemann surfaces.
\end{abstract}


%%%%%%%%%%%%%%%%%%%%%%%%%%%%%%%%%%%%%
%									%
%			\maketitle		 		%
%									%
%%%%%%%%%%%%%%%%%%%%%%%%%%%%%%%%%%%%%
\maketitle

%%%%%%%%%%%%%%%%%%%%%%%%%%%%%%%%%%%%%
%									%
%		Table Of Contents			%
%									%
%%%%%%%%%%%%%%%%%%%%%%%%%%%%%%%%%%%%%

\setcounter{tocdepth}{1}
\tableofcontents

%%%%%%%%%%%%%%%%%%%%%%%%%%%%%%%%%%%%%
%									%
%			Introduction			%
%									%
%%%%%%%%%%%%%%%%%%%%%%%%%%%%%%%%%%%%%

\section{Introduction}
\setcounter{equation}{0}

In the present paper we extend the theory of pairs of commuting non-selfadjoint operators in Hilbert spaces as presented in
\cite{KLMV} to the Pontryagin space setting.
To set the framework and provide
motivation we first discuss the case of a single operator in a Hilbert space and
we recall that models of linear operators in Hilbert space and reproducing kernel Hilbert spaces of analytic functions,
in particular of the type introduced by de Branges \cite{dbbook} and
de Branges and Rovnyak \cite{dbr2}, are closely related topics.
Typically (but not always), $A$ will be unitary equivalent to the backward-shift operator
\[
R_0 \, f (z) = \frac{f(z) - f(0)}{z}
\]
in some reproducing kernel Hilbert space of analytic functions.\smallskip

In the case of a non-selfadjoint (or non-unitary) operator,
one associates in a natural way a Schur function,
i.e. a function analytic and contractive
(possibly with respect to an indefinite metric) in the open unit disk or the open half-plane and such that
the kernel
\begin{equation}
\label{123456}
\frac{J-S(z)JS(w)^*}{1-z\overline{w}},
\end{equation}
or
\begin{equation}
\frac{-J+S(z)JS(w)^*}{-i(z-\overline{w})},
\label{12345}
\end{equation}
is positive definite in the open unit disk (respectively in the open upper-half-plane) from which is possibly removed a zero set.
Here $J$ is a signature matrix (i.e. both selfadjoint and unitary).
Note that in the work of de Branges and other places the opposite of the kernel \eqref{12345} is considered,
here we follow the convention in \cite{brod}.
We remark that analyticity of $S$ and positivity of the corresponding kernel in an open subset is equivalent to the $J$-contractivity of $S$ in that subset. In particular,
if $J=1$ and $S$ is analytic and contractive, the kernels above are positive definite in the open unit disk or the upper half-plane
(as we will see below, this is not necessarily true when one considers the Riemann surface case).
\smallskip

One can define the characteristic function in various ways \cite{brod,nf}.
Here for a contraction, say $T$, we have in mind the transfer operator of the associated Julia operator, i.e.
\begin{equation}
\label{3456}
S(z)=
H + z G (I-z T)^{-1}F
\end{equation}
where the matrix operator
\[
\begin{pmatrix}
T&F\\ G&H
\end{pmatrix}
\]
is unitary and is called (when uniquely defined from $T$) a Julia operator associated to $T$.\\

Brodski{\u\i} and Liv{\v{s}}ic \cite{MR20:7221}, Liv{\v{s}}ic \cite{l1} and Brodski{\u\i} \cite{MR48:904} 
presented a system theoretic interpretation to Schur functions.
The function $S$ in \eqref{3456} is the transfer function of a corresponding dissipative $1D$ linear system
\begin{align*}
   i\frac{d}{dt}f(t) +A f(t) &=\Phi^* \sigma u(t),\\
    y(t)&=u(t)-i \Phi f(t)
.
\end{align*}
Here $f(t)$ is the state signal with values in a Hilbert space $H$, $u(t)$ and $y(t)$ are the input and output signals, respectively, 
with values in the coefficient (Hilbert) space $E$, $\Phi$ is a bounded operator from $H$ to $E$ and $\sigma \in \mathcal B (E)$.
\smallskip

In the half-plane case, the associated kernel is \eqref{12345}
and the corresponding characteristic operator is of the form
\[
S(z) = I - i \Phi (A- zI)^{-1} \Phi^* J
,
\]
satisfying
\[
\frac{A - A^*}{i} = \Phi^* J \Phi,
\]
and the following identity holds
\[
-J + S(z) J S(w)^*
=
i (z- \overline{w}) \Phi (A - Iz)^{-1}(A^* - I\overline{w})^{-1} \Phi^*
.
\]
Single operator models and the corresponding system theory have been generalized to the Pontryagin case,
see \cite{MR2002b:47144,adrs,ad3}, \cite{MR2269246,MR2003e:47022} and \cite{MR1664343,MR1671482,MR1991668}
and others.
\\

The case of operators in Pontryagin spaces can be motivated by the theory of
ordinary differential operators with certain boundary conditions \cite{MR2841147,MR3052315},
but is also a topic of interest in its own right.
Schur functions are replaced by generalized Schur functions.
Recall that when $J=I$, a meromorphic function $S(z)$ such that \eqref{123456}
has a finite number of negative squares if and only if $S(z)$
can be written as $S(z)=S_0(z) B_0^{-1}(z)$, where $S_0(z)$ is Schur function
and $B_0(z)$ is a finite Blaschke product; see \cite{kl1}.
This will not be true in general in the Riemann surface setting.
See the discussion below.
\smallskip

When going from a single non-selfadjoint operator to a pair of commuting non-selfadjoint operators in a Hilbert space,
the plane (or more precisely the Riemann sphere) is replaced by a real compact Riemann surface
(a compact Riemann surface equipped with an anti-holomorphic involution); see \cite{KLMV}.
One still has the notion of Schur function.
In the scalar case, the counterpart of a reproducing kernel of the form \eqref{12345}, is the following expressions (see \cite{vinnikov4})
\begin{equation}
\label{rsKernel}
\frac{\vartheta[\widetilde{\zeta}](\overline{q}-p)}
{i \vartheta[\widetilde{\zeta}](0)E(p,\overline{q})}
-
T(p)
\frac{\vartheta[\zeta](\overline{q}-p)}{i \vartheta[\zeta](0)E(p,\overline{q})}
\overline{T(q)}
,
\end{equation}
where $\vartheta[\zeta]$ is the theta function with characteristic $\zeta$
and $E(p,q)$ is the prime form; we refer to \cite{mumford1,mumford2} and \cite{fay1} for more details.
These kernels were studied in the finite dimensional case in \cite{av3}.
As shown in \cite{av3}, even in the scalar and dividing case,
where $X \backslash X_{\mathbb R}$ contains two connected components, Pontryagin spaces appear naturally.
Even if an analytic bundle mapping on $X_+$ is contractive, the associated kernel is not necessarily positive.
In opposition with the complex setting,
there exists $T$ such that the linear span of the sections of the form \eqref{rsKernel} where $q \in X \backslash X_{\mathbb R}$, 
is finite dimensional, but $T$ is not a quotient
$T = T_1 T_2^{-1}$ of two finite Blaschke products. See \cite{av3} for an example.
\smallskip

Also in the $2D$ case there is a system theory interpretation.
The characteristic function of an overdetermined $2D$ system becomes a mapping between certain vector bundles associated to a curve,
called the discriminant curve. See for instance, \cite{livsic3,livsic1,KLMV}.
In the upcoming pages we develop the corresponding indefinite theory and consider overdetermined $2D$ systems in the non-positive case.
In particular, we develop the vessel theory in Pontryagin spaces.

%%%%%%%%%%%%%%%%%%%%%%%%%%%%%%%%%%%%%
%									%
%	Colligations					%
%									%
%%%%%%%%%%%%%%%%%%%%%%%%%%%%%%%%%%%%%
\subsection{Structure of the paper}
This paper consists of eight sections besides the introduction and we now sketch briefly its contents.
\smallskip

In Section 2, we present a short survey of indefinite inner product spaces, and in particular, Pontryagin spaces.
To go from the Hilbert space setting
to the Pontryagin space case, we make use of a number of important facts which
still hold in the Pontryagin space case;
see Theorems \ref{IKLthm1}-\ref{IKLthm3} and Theorems \ref{IKLthm4}-\ref{IKLthm7}.
\smallskip

The corresponding single-operator colligation theory is developed and presented in Section 3.
In particular, we give a realization theorem of function $S$ such that the corresponding kernel \eqref{12345} has finite number of negative squares
and analytic at infinity.
In Section 4, we define and study the vessel theory over Pontryagin spaces.
Section 5 is dedicated to the corresponding $2D$ systems.
The coupling and the decompositions of $2D$ systems, or equivalently, vessels over Pontryagin space are studied in Section 6.
\smallskip

Section 7, Section 8 and Section 9 contain the main results of the paper.
The various characteristic functions and the associated realization theorems are presented in Section 7.
In Section 8, we deduce a de Branges spaces characterization theorem in the non-positive case.
Finally, in Section 9, we review the indefinite Hardy spaces of sections defined on finite bordered Riemann surfaces (as developed in \cite{av2}) and use 
the de Branges theory to present a characterization theorem (Beurling's Theorem) 
of the invariant subspaces under the multiplication operators in the indefinite setting.

%%%%%%%%%%%%%%%%%%%%%%%%%%%%%%%%%%%%
%									%
%	Colligations					%
%									%
%%%%%%%%%%%%%%%%%%%%%%%%%%%%%%%%%%%%
\begin{Rk}
Theorems \ref{IKLthm1}-\ref{IKLthm3} and Theorems \ref{IKLthm4}-\ref{IKLthm7}
(which as already mentioned play an important role in the arguments)
do not hold in the setting of Krein spaces, which requires supplementary hypothesis.
The counterpart theory for Krein spaces will be elaborated in a future publication.
\end{Rk}
\begin{Rk}
In this paper we focus on the case of a pair of commuting operators. The generalization to the case of $n$-tuple of commuting operators will be also considered elsewhere.
\end{Rk}

%%%%%%%%%%%%%%%%%%%%%%%%%%%%%%%%%%%%
%									%
%	Colligations					%
%									%
%%%%%%%%%%%%%%%%%%%%%%%%%%%%%%%%%%%%
\subsection{Some notations:}
Pontryagin spaces are denoted by $\mathcal P$ with the finite negative index $\kappa$ or $\kappa _{\mathcal P}$.
We use the symbol $\mathcal K$ to denote a Krein space.
The inner product over a Hilbert space is denoted by $\innerProductReg{\cdot}{\cdot}$,
while $\innerProductIndefiniteReg{\cdot}{\cdot}$ is the indefinite inner product.
\smallskip

The adjoint operator over Pontryagin space is denoted by $A^{[*]}$ and by $A^*$ in the Hilbert space case.
The complete, joint and normalized joint characteristic functions of a vessel are denoted by $W(\xi_1,\xi_2,z)$, $S(\lambda)$ and $T(x)$, respectively.

%%%%%%%%%%%%%%%%%%%%%%%%%%%%%%%%%%%%%
%									%
%	Preliminaries					%
%									%
%%%%%%%%%%%%%%%%%%%%%%%%%%%%%%%%%%%%%

\section{Preliminaries}
\setcounter{equation}{0}

%%%%%%%%%%%%%%%%%%%%%%%%%%%%%%%%%%%%%
%									%
%	Pontryagin Spaces				%
%									%
%%%%%%%%%%%%%%%%%%%%%%%%%%%%%%%%%%%%%

\subsection{Indefinite inner product spaces}
This section is dedicated to a brief overview of indefinite inner product spaces.
We focus, in particular, on Pontryagin spaces and their fundamental properties.
For general background on indefinite inner product spaces, especially Pontryagin and Krein Spaces, we refer the readers to \cite{azih,bognar,ikl}.
\smallskip

Let us consider a vector space $\mathcal V$ over the complex numbers, endowed with an Hermitian form
$ \innerProductIndefiniteReg{\cdot} {\cdot}$.
Given two linear subspaces $\mathcal V _1$ and $\mathcal V _2$  such that $\mathcal V _1 \cap \mathcal V _1 = \{0\}$, we denote their direct sum by $\mathcal V _1 \oplus \mathcal V _1$.
Two elements $v,w \in \mathcal V$ are called orthogonal with respect to this form
if $\innerProductIndefiniteReg{v}{w} = 0$. Two linear subspaces $\mathcal V _1$ and $\mathcal V _2$
are orthogonal if every element of the first subspace is orthogonal to every element of the second. We use the notation $\mathcal V _1 [+] \mathcal V _2$ to denote the orthogonal sum.
If moreover $\mathcal V _1 \cap \mathcal V _1 = \{0\}$, the sum (which is then also direct) is denoted
by
\[
\mathcal V _1 [ \oplus ] \mathcal V _1.
\]
\begin{definition}
The space $\mathcal V$ endowed with the Hermitian form $\innerProductIndefiniteReg {\cdot} {\cdot}$ is called a Krein
space if there exist two subspaces $\mathcal V _1$, $\mathcal V _2$ such that $\mathcal V$ can be written as an orthogonal
and direct sum
\begin{equation}
\label{kreinDeco}
\mathcal V _+ [ \oplus ] \mathcal V _-,
\end{equation}
where ($\mathcal V _+$; $\innerProductIndefiniteReg {\cdot} {\cdot }$) and ($\mathcal V _-$; $-\innerProductIndefiniteReg {\cdot} {\cdot }$) are both Hilbert spaces. The space $\mathcal V$ is called a Pontryagin
space if $\mathcal V _-$ is finite dimensional. The dimension of the vector space $\mathcal V _-$ is called the (negative) index of the Pontryagin space.
\end{definition}

The decomposition in \eqref{kreinDeco}, called a fundamental decomposition, will not be unique, unless
one of its component equals the subspace $\{ 0 \}$. The signature operator is the operator $J_{\mathcal V} : \mathcal V \rightarrow \mathcal V$ such that
\[
J_{\mathcal V} (v) = J_{\mathcal V} (v_+ + v_-) = v_+ - v_- \qquad ,  \qquad \forall v_{\pm} \in \mathcal V _{\pm}
.
\]
A Krein space $\mathcal V$ becomes a Hilbert space when endowed with the inner product:
\begin{equation}
\label{innerProd}
\innerProductReg{v}{w} = \innerProductIndefiniteReg{v_+}{w_+} - \innerProductIndefiniteReg{v_-}{w_-},
\end{equation}

where $v = v_+ + v_-$ and $w = w_+ + w_-$ belong to $\mathcal V$ and $v_\pm,w_\pm \in \mathcal V _\pm$. The inner product
\eqref{innerProd} depends on the given decomposition, but all the resulting norms are equivalent, and
defines the topology of the Krein space (see \cite[p. 102]{bognar}). The notion of continuity and
of convergence are with respect to this topology. Given two Krein spaces $\mathcal V$, $\mathcal W$ we denote
by $\mathcal B (\mathcal V,\mathcal W )$ the set of continuous linear operators from $\mathcal V$ to $\mathcal W$ and when $\mathcal W = \mathcal V$
we will use the symbol $\mathcal B(\mathcal V )$. By $I_{\mathcal V}$ we denote the identity operator on $\mathcal V$.

\begin{Rk}
In the finite dimensional case, let $\mathcal V = \mathbb C ^m$ and let $J \in \mathbb C^{m \times m}$ be a signature matrix. Then $ \mathcal V$ endowed with the Hermitian form
\[
\innerProductIndefiniteReg{u}{v} = v^*Ju, \qquad u,v \in \mathbb C ^m,
\]
is a Pontryagin space and a fundamental symmetry is given by the map $u \rightarrow Ju$.
\end{Rk}

The adjoint of an operator can be defined with respect to the Hilbert space inner
product or with respect to the Krein space inner product. More precisely, let $(\mathcal V , \innerProductReg{\cdot}{\cdot}_{\mathcal V})$ and
$(\mathcal W , \innerProductReg{\cdot }{ \cdot}_{\mathcal W})$ be Hilbert spaces. Given $A \in \mathcal B( \mathcal V , \mathcal W )$, 
its adjoint is the unique operator $A^{[*]} \in \mathcal B( \mathcal W , \mathcal V )$ such that
\[
\innerProductReg{Af }{ g}_{\mathcal W}
=
\innerProductReg{f } {A^{*} g}_{\mathcal V},
\quad
{\rm where} \quad f \in \mathcal V
\quad
{\rm and} \quad g \in \mathcal W.
\]

\begin{definition}
Let $(\mathcal V , \innerProductIndefiniteReg{\cdot}{ \cdot}_{\mathcal V})$ and $(\mathcal W , \innerProductIndefiniteReg{\cdot}{\cdot}_{\mathcal W})$ be two Krein spaces. Given $A \in \mathcal B(\mathcal V, \mathcal W )$ its adjoint is the unique operator $A^{[*]} \in \mathcal B(\mathcal W ,\mathcal V )$ such that
\[
\innerProductIndefiniteReg{Af}{ g}_{\mathcal W}
=
\innerProductIndefiniteReg{f}{ A^{[*]}g}_{\mathcal V} ,\qquad f \in \mathcal V \quad g \in \mathcal W.
\]
\end{definition}

The Krein spaces $\mathcal V$ and $\mathcal W$ become Hilbert spaces when endowed with the Hermitian forms
\begin{equation}
\label{kreinHilbertProduct}
\innerProductReg{ f_1}{ f_2 }_{\mathcal V} = \innerProductIndefiniteReg{ f_1}{J_{\mathcal V} f_2}_{\mathcal V}
\qquad {\text and } \qquad
\innerProductReg{g_1}{g_2}_{\mathcal W} = \innerProductIndefiniteReg{g_1}{J_{\mathcal W} g_2}_{\mathcal W},
\end{equation}
where $f_1, f_2 \in \mathcal V$ and  $g_1,g_2 \in \mathcal W$.
\smallskip

The Hilbert space adjoint operator $A^*$ (with respect to the inner products \eqref{kreinHilbertProduct}) and the Krein space adjoint operator $A^{[*]}$ are related by
\begin{equation}
\label{adjoint2adjoint}
A^* = J_{\mathcal V} A^{[*]} J_{\mathcal W}
,
\end{equation}
where $J_{\mathcal V}$ and $J_{\mathcal W}$ are associated to some fundamental decompositions of $\mathcal  V$ and $\mathcal  W$ .

\begin{definition}
Let $\mathcal V$ and $\mathcal W$ be two Krein spaces. 
An operator $A \in \mathcal B(\mathcal V , \mathcal W )$ is said to be:
\begin{enumerate}
\item isometric if $A^{[*]}A = I_{\mathcal V}$.
\item coisometric if $AA^{[*]} = I_{\mathcal W}$.
\item unitary if it is isometric and coisometric.
\item a contraction if $\innerProductIndefiniteReg{Af}{ Af }_{\mathcal W}  \leq \innerProductIndefiniteReg{ f }{ f }_{\mathcal V} ,$ for all  $f \in \mathcal V$.
\end{enumerate}
\end{definition}
Here we mention several important results which are used in this paper.
We note that some of these statements hold only in the Pontryagin space case, but fail in Krein space setting 
(we refer to \cite{ikl} for the proofs).
\begin{theorem}
\label{IKLthm1}
The adjoint of a contraction between Pontryagin spaces of same index
is also a contraction.
\end{theorem}
\begin{theorem}[{\cite[Theorem 1.4.1]{adrs}}]
\label{IKLthm3}
A densely defined contractive operator (resp. an isometry) between Pontryagin spaces of
same index extends to a continuous everywhere defined contraction (isometry).
\end{theorem}

We present below some more related and relevant results.

\begin{theorem} [{\cite[Theorem 2.4]{ikl}}]
\label{IKL24}
\label{IKLthm4}
Let $\mathcal P$ be a Pontryagin space of $\kappa$ negative index, and let $\mathcal L$ be an arbitrary dense subset in $\mathcal P$.
\begin{enumerate}
\item
The sequence $ \left( x_n\right)_{n=1}^{\infty}$ where $x_n \in \mathcal P$
converges to some $x_0 \in \mathcal P$ with respect to an arbitrary norm if and only if for all $y \in \mathcal L$
\begin{equation}
\label{csConverge}
    \innerProductIndefiniteReg{x_n}{y} \xrightarrow{n \rightarrow \infty} \innerProductIndefiniteReg{x_0}{y}
\end{equation}
and
\[
\innerProductIndefiniteReg{x_n}{x_n}  \xrightarrow{n \rightarrow \infty} \innerProductIndefiniteReg{x_0}{x_0}.
\]
\item
The sequence $ \left( x_n\right)_{n=1}^{\infty}$ where $x_n \in \mathcal P$
is a Cauchy sequence with respect to an arbitrary norm if and only if
\[
\innerProductIndefiniteReg{x_n - x_m}{x_n - x_m} \xrightarrow{n \rightarrow \infty} 0.
\]
and $\left( \innerProductIndefiniteReg{x_n}{y} \right)_{n=1}^{\infty}$ is a Cauchy sequence for all elements $y$ in a dense subset of $\mathcal L$.
\end{enumerate}
\end{theorem}
The fact that $\mathcal L$ may be chosen dense (but $\neq \mathcal P$)
is of key importance in the reproducing kernel Pontryagin space setting.
One can then take $\mathcal L$ to be the linear span of the kernel and then
\eqref{csConverge} becomes pointwise convergence.

We conclude this section with three results used in the sequel, but first we recall that
a subspace $\mathcal L$ is called ortho-complemented if it admits an orthogonal complementary subspace.
\begin{theorem}[{\cite[Theorem 2.2  P.186]{bognar}}]
\label{bognar22}
\label{IKLthm5}
Every closed, non--degenerate subspace $\mathcal L$ of a Pontryagin space $\mathcal P$ with negative index $\kappa$, is ortho-complemented.
\end{theorem}

\begin{theorem}[{\cite[Corollary 2.3  P.186]{bognar}}]
\label{bognar23}
\label{IKLthm6}
Every closed, non--degenerate subspace of a Pontryagin space $\mathcal P$ with negative index $\kappa$, is a Pontryagin space $\mathcal P_1$ with negative rank $\kappa_1 \leq \kappa$.
\end{theorem}

To put the previous two theorems in perspective,
we mention the following result which will not be used in sequel.

\begin{theorem}[{\cite[Theorem 3.4 P.104]{bognar}}]
\label{IKLthm7}
A subspace $\mathcal L$ of the Krein space $\mathcal K$
is closed and a Krein space if and only if $\mathcal L$ is ortho-complemented.
\end{theorem}

\subsection{Reproducing kernel Pontryagin space}
\begin{definition}
Let $\Omega$ be some set and let $\mathcal K$ be a Krein space. The $\mathcal B( \mathcal K )$-valued function $K(z,w)$ defined on $\Omega \times \Omega$ is said to have $\kappa$ negative squares if it is Hermitian, that is,
\[
K(z,w) = K(w, z)^{[*]}, \qquad \forall z,w \in \Omega
\]
and if for every choice of $N \in \mathbb N, c_1,...,c_N \in \mathcal K$ and $w_1,...,w_N \in \Omega$ the $N \times N$ Hermitian
matrix with $(u,v)$-entry equal to
\[
\innerProductIndefiniteReg{K(w_u,w_v)c_u}{c_v}_{\mathcal K}
\]
has at most $\kappa$ strictly negative eigenvalues, and exactly $\kappa$ such eigenvalues for some choice of $N,c_1,...,c_N$ and $w_1,...,w_N$.
\end{definition}

\begin{Rk}
Here as in \cite{brod} and \cite{kl1} 
the matrix-valued kernels $K(z,w)$ are analytic from the right side in $z$
and from the left side in $\overline{w}$.
When the opposite convention is taken,
one replaces
$\innerProductIndefiniteReg{K(w_u,w_v)c_u}{c_v}_{\mathcal K}$
by
    \[\innerProductIndefiniteReg{K(w_u,w_v)c_v}{c_u}_{\mathcal K}.\]
\end{Rk}

The function $K(z,w)$ is called positive definite if $\kappa = 0$
(here we say positive rather than positive definite).
We refer to the function $K(\cdot,\cdot)$ as kernel. 
An important related notion is the reproducing kernel Pontryagin space associated to $K(z,w)$.

\begin{definition}
Let $\Omega$ be some set, let $\mathcal K$ be a Krein space and let $\mathcal P$ be a Pontryagin space of $\mathcal K$-valued functions defined on $\Omega$.
Then $\mathcal P$ is called a reproducing kernel Pontryagin space if there exists a $\mathcal B(\mathcal K )$-valued function $K(z,w)$ with the following two properties: For every $c \in \mathcal K , w \in \Omega$ and $F \in \mathcal P$,
\begin{enumerate}
\item The function $z \rightarrow K(z,w)c$ belongs to $\mathcal P$.
\item It holds that
\[
\innerProductIndefiniteReg{F}{K(\cdot,w)c}_{\mathcal P} = \innerProductIndefiniteReg{F(w)}{c}_{\mathcal K}.
\]
\end{enumerate}
\end{definition}
The function $K(z,w)$ is the reproducing kernel associated to $\mathcal P$.
It is Hermitian and by the Riesz' representation theorem it is uniquely defined.
\begin{definition}
\label{defRank}
The rank of a Hermitian kernel $K(z,w)$ is the supremum
over the ranks of the matrices $\left(K(z_i,z_j)\right)_{i,j \in I}$
where $I$ is a finite subset $I = \{i_1,...,i_n\}\subset \Omega$.
The rank of a kernel, by definition, is either finite or infinite.
\end{definition}
The following theorems relate the two definitions above, and originates with the work of Aronszajn \cite{aron}
in the case of positive definite kernels and Schwartz \cite{Sch95} and Sorjonen \cite{sorjonen} in the case of negative squares.
\begin{theorem}
\label{pontKernelDecomB}
Let $\Omega$ be some set and let $\mathcal K$ be a Krein space.
There is a one-to-one correspondence between reproducing kernel Pontryagin spaces of negative index $\kappa$ of $\mathcal K$-valued functions defined on $\Omega$
and $\mathcal B(\mathcal K )$-valued functions having $\kappa$ negative squares in $\Omega$.
\end{theorem}
\begin{theorem}
\label{pontKernelDecomA}
A Hermitian kernel $K(z,w)$ on a set $\Omega$ with values in ${\mathcal B}({\mathcal K})$ has $\kappa$ negative squares if and only if
there exists a Pontryagin space ${\mathcal P}$ of negative index $\kappa$ and a function $L \colon \Omega \to {\mathcal B}({\mathcal P},{\mathcal K})$
with $$ K(z,w) = L(z) L(w)^{[*]}$$ and $\bigvee_{w \in \Omega} \operatorname{Im} L(w)^{[*]} = {\mathcal P}$.
\end{theorem}
We note that sufficiency follows from Theorem \ref{pontKernelDecomB} by taking $L(z)$ to be the evaluation at $z$.
\begin{theorem}
The function $K(z,w)$ has $\kappa$ negative squares in $\Omega$ if and only if it can be written as $K(z,w) = K_+(z,w)-K_-(z,w)$,
where both $K_+(z,w)$ and $K_-(z,w)$ are positive definite in $\Omega$ and moreover $K_-$ has rank $\kappa$ (see Definition \ref{defRank}).
\end{theorem}

%%%%%%%%%%%%%%%%%%%%%%%%%%%%%%%%%%%%%
%									%
%									%
%									%
%%%%%%%%%%%%%%%%%%%%%%%%%%%%%%%%%%%%%

\section{Single-operator colligation over Pontyagin space}
\setcounter{equation}{0}

A characterization of reproducing kernel Pontryagin spaces
with reproducing kernel of the form \eqref{123456}
has been considered in \cite{adrs} using the theory of linear relations in Pontryagin space.
Here we study kernels of the form \eqref{12345} using different methods.
%%%%%%%%%%%%%%%%%%%%%%%%%%%%%%%%%%%%%
%									%
%	Realization						%
%									%
%%%%%%%%%%%%%%%%%%%%%%%%%%%%%%%%%%%%%

\subsection{Realization of $J$-unitary function analytic at infinity}

We begin this section recalling the following lemma.
Let $S$ be a function defined on a set $\Omega$ with values in a
Hilbert space $E$.
We denote by $\mathcal P _K(S)$ the reproducing kernel Pontryagin space with kernel of the form
\[
K_w (z)
=
\frac{-J+S(w)^{*}JS(z)}{-i(z- \overline{w})}
,
\]
where $J$ is a selfadjoint and unitary operator.
\begin{lem}[{\cite[Theorem 6.9]{ad3}}]
Let $S$ be an analytic function on an open set $\Omega$ with values in a
Hilbert space $E$ such that $S(\overline{z})^* J S(z) = J$ whenever $z,\overline{z} \in \Omega$.
Then
$\mathcal P_K (S)$ is $R_\alpha$-invariant where $\alpha \in \Omega$ .
\end{lem}
Realizations for elements in $\mathcal P_K(S)$ centered at a finite point
can be found in \cite{ad3} and \cite{alpay2017beurling}.
Here we consider realizations as $\alpha$ tends to infinity and provide direct argument.
\begin{theorem}
\label{th1}
Let $S$ be a function with values in a
Hilbert space $E$ analytic a neighborhood of infinity, say $\Omega$, such that $S(\infty)=I$ and
$S(\overline{z})^* J S(z) = J$ whenever $z,\overline{z} \in \Omega$.
Furthermore, let us assume that the kernel $K(z,w)$ has $\kappa$ negative squares.
Then $\mathcal P_K (S)$ does not contain nonzero constants
and $S (z)$ has a realization of the form
\begin{equation}
\label{thetaReal}
S(z) = I + i J C (zI-A^\prime)^{-1} C^{[*]},
\end{equation}
    for some operator $A^\prime \in  \mathcal B(H)$ and $C \in \mathcal B (H,E)$ satisfying
    \begin{equation}
        \label{eqCollForS}
\frac{1}{i}
\left(
A^\prime - A ^{^\prime [*]}
\right)
= C^{[*]} J C.
    \end{equation}
\end{theorem}
\begin{Rk}
Using the Potapov-Ginzburg transform as in \cite[Theorem 4.3.5]{adrs},
one sees that $S$ has an extension of bounded type to the upper half-plane ($\mathbb{C}_+$).
Thus $S(z)$ admits,
almost everywhere on the real line,
non tangential limits
(which we denote by $S(x)$).
Then it makes sense  to consider an additional property,
namely that $S$ is $J$-unitary on the real line:
\[
S(x)^{*} J S(x)  = J \, \, \, \, x \in \mathbb{R}.
\]
\end{Rk}
\begin{pf}[of Theorem \ref{th1}]
    We use the Potapov-Ginzburg transform to map $S$ to $\Sigma$, a meromorphic function defined on a Hilbert space.
    The transform is given explicitly by
\[
\Sigma(\lambda) = \left(P S(\lambda) +  Q \right) \left( P + Q S(\lambda) \right)^{-1}
.
\]
Here we define the following pair of selfadjoint projections
\begin{equation}
\label{eqPq}
P \defEq \frac{I+J}{2} \quad {\rm and} \quad Q \defEq \frac{I-J}{2},
\end{equation}
We note that since $S(\infty)=I$, then it follows that $P +Q S(\lambda) $ is invertible in the neighborhood of infinity.
Furthermore, obviously $\Sigma(\infty)=I$ and the kernel $K_w(z)$ can be written as
\begin{equation}
\label{eqKernelSKernelOmega}
\frac{J-S(w)^* J S(z)}{-i(z- \overline{w})}
=
\left(P + Q S(w)\right)^*
\frac{I-\Sigma(w)^* \Sigma(z)}{-i(z- \overline{w})}
\left(P + Q S(z) \right) 
.
\end{equation}
At this stage, we divide the proof into a number of steps.
\smallskip

    \begin{Step}
        We prove \eqref{eqSigmaReal}  for a realization of $\Sigma$.
    \end{Step}
It follows from \eqref{eqKernelSKernelOmega} that $\Sigma(z)$ is a generalized Schur function,
    that is, the kernel
    \[
        K(z,w) = \frac{I- \Sigma(w)^* \Sigma(z) }{-i(z-\overline{w})},
    \]
    has a finite number of negative squares.
Hence, by the Krein-Langer factorization theorem, we have the following decomposition
$$\Sigma(z) = B_0 ^{-1}(z) \Sigma_0(z).$$
Here $B_0(z)$ is a Blaschke product of degree $\kappa$ and $\Sigma_0(z)$ belongs to the Schur class (see \cite{kl1} and \cite{adrs}).
The matrix function $\Sigma_0(z)$, as a Schur function, admits a realization (see for instance \cite{brod}) 
\[
\Sigma_0(z) = D_1 + C_1 (zI - A_1)^{-1} B_1
\]
which satisfy 
    \begin{equation}
        \label{eqRealSigma0}
        B_1 = i C_1, \qquad A_1 - A_1^* = i C_1^* C_1,
    \end{equation}
    and $D_1=I $ since $\lim_{z \rightarrow \infty} \Sigma_0(z) = I$.
The finite Blaschke product $B_0(z)$ also has a realization
\[
B_0(z) = D_3 + C_3 (zI - A_3)^{-1}B_3
,
\]
satisfying 
\begin{equation}
    \label{eqRealB}
    B_3 = i C_3, \qquad A_3 - A_3^* = i C_3^* C_3,
\end{equation}
    and $D_3=I $ since $\lim_{z \rightarrow \infty} B_0(z) = I$.
\smallskip

Then, a realization of $B_0 ^{-1}(z)$ is given by the block matrix
\begin{equation}
    \label{eqRealationB}
\begin{pmatrix}
  A_2 & B_2
  \\
  C_2 & D_2
\end{pmatrix}
=
\begin{pmatrix}
  A_3-B_3 D_3^{-1} C_3 & B_3 D_3^{-1}
  \\
  -D_3^{-1} C_3 & D_3^{-1}
\end{pmatrix}
=
\begin{pmatrix}
  A_3-B_3 C_3 &B_3
  \\
  - C_3 & I
\end{pmatrix}
\end{equation}
and hence, substituting \eqref{eqRealationB} in \eqref{eqRealB}, we have
\begin{equation}
    \label{eqRealBinverse}
    B_2^{ *} = -iC_2, \quad {\rm  and} \quad A_2 - A_2^* = -i C_2^* C_2.
\end{equation}
Moreover, by the theorem of multiplication of two realizations (see for instance \cite[Section 2.5]{MR2663312}) we have
\begin{equation}
\label{eqSigmaReal1}
\Sigma(z) = \Sigma(\infty) + C(z I - A )^{-1}B,
\end{equation}
where
\[
    A = 
\begin{pmatrix}
    A_1 & B_1 \, C_2  \\
    0   & A_2
\end{pmatrix}
\]
and
\begin{equation}
    \label{eqBandC}
B =
\begin{pmatrix}
  B_1 D_1
  \\
    B_2
\end{pmatrix}
=
\begin{pmatrix}
  B_1
  \\
    B_2
\end{pmatrix}
\qquad
{\rm and}
\qquad
C =
\begin{pmatrix}
  C_1
  &
    D_1 C_2
\end{pmatrix}
=
\begin{pmatrix}
  C_1
  &
  C_2
\end{pmatrix}
.
\end{equation}
The matrices $B$ and $C$ are related by
\[
C^*
=
\begin{pmatrix}
  C_1^*
  \\
  C_2^*
\end{pmatrix}
=
\begin{pmatrix}
  iB_1
  \\
  -iB_2
\end{pmatrix}
=
i
\begin{pmatrix}
 I & 0  \\
 0 & -I \\
\end{pmatrix}
B
.
\]
Hence the adjoint of the operator $C$ relative to the indefinite inner product
defined by $\begin{psmallmatrix}I&0\\0&-I \end{psmallmatrix}$, is given by
\begin{equation}
\label{eqCB}
C^{[*]} = iB,
\end{equation}
and \eqref{eqSigmaReal1} becomes
\begin{equation}
\label{eqSigmaReal}
\Sigma(z) = I + i B^{[*]}(z I - A )^{-1}B.
\end{equation}
Finally, we have the following 
\begin{align*}
    A - A^{[*]} & = A - J A^* J =
\begin{pmatrix}
    A_1 & B_1 \, C_2  \\
    0   & A_2
\end{pmatrix}
-
\begin{pmatrix}
    A_1^* & 0  \\
    - C_2^* B_1^*   & A_2^*
\end{pmatrix}
.
\end{align*}
Using \eqref{eqRealSigma0}, \eqref{eqRealBinverse} and \eqref{eqBandC}, we have
\begin{align}
    \nonumber
    A - A^{[*]} 
    & =
\begin{pmatrix}
    A_1-A_1^* & B_1 \, C_2  \\
    C_2^* B_1^*   & A_2 - A_2^*
\end{pmatrix}
=
\begin{pmatrix}
    iC_1^*C_1 & B_1 \, C_2  \\
    C_2^* B_1^*   & -iC_2^* C_2
\end{pmatrix}
\\    & =
\begin{pmatrix}
    iC_1^*C_1           & iC_1^* \, C_2  \\
    -i C_2^* C_1      & -iC_2^* C_2
\end{pmatrix}
    =iJ C^* C = iC^{[*]} C,
    \label{eqRealSigma}
\end{align}
and the colligation condition for $\Sigma$ follows.
\begin{Step}
    We prove the realization formula \eqref{eq1234} for $S$.
\end{Step}
We rewrite $S$ in term of $\Sigma$ (the inverse Potapov-Ginzburg transform)
\begin{align}
S(z)
= &
\label{eqSinQPSIGMA}
    \left(Q  + P \Sigma(z)  \right)  \left(P + Q \Sigma(z) \right) ^{-1}
\\
= &
\nonumber
\left(Q  + P  B_0 ^{-1}(z) \Sigma_0(z)  \right)  \left(P+ Q B_0 ^{-1}(z) \Sigma_0(z) \right)
.
\end{align}
We substitute \eqref{eqPq} and \eqref{eqSigmaReal} in \eqref{eqSinQPSIGMA} and have the following
\begin{align*}
S(z)
= &
    \left(\frac{I-J}{2}  +  \frac{I+J}{2}\Sigma(z) \right) \left(  \frac{I+J}{2} + \frac{I-J}{2}\Sigma(z) \right)^{-1}
\\
= &
\left(\frac{I-J}{2}  + \frac{I+J}{2} \left( I + C(z I - A )^{-1}B \right)  \right)
\times
\\
&
    \left( \frac{I+J}{2} + \frac{I-J}{2} \left( I + C(z I - A )^{-1}B \right)  \right)^{-1}
\\
= &
    \left(2 I + (I-J+2J)C(z I - A )^{-1}B  \right)
\times
\\ &
    \left(2 I + (I-J)C(z I - A )^{-1}B  \right)^{-1}
.
\end{align*}
Using the equality $F^{-1} B (I+DF^{-1}B)^{-1} = (F+BD)^{-1} B $, we then have
\begin{align*}
S(z)
= &
I +
 J C(z I - A )^{-1}B
    \left( I + \frac{I-J}{2} C(z I - A )^{-1}B   \right)^{-1}
    \\
= &
I +
J C 
\left( (z I - A ) + B  \frac{I-J}{2} C \right)^{-1}B
.
\end{align*}
Thus $S(z)$ has a realization of the form
\begin{equation}
    \label{eq1234}
S(z)
=
I +
J C \left( z I - A^{\prime}\right)^{-1}B,
\end{equation}
where $A^{\prime} = A  - B Q C$. Furthermore, recalling that \eqref{eqCB} holds, \eqref{thetaReal} follows.
Note that \eqref{eqCollForS} 
is a consequence of the colligation condition associated to the realization of $\Sigma$ (i.e.  $A-A^{[*]} = iC^{[*]} C$, see \eqref{eqRealSigma}), 
    $iB=C^{[*]}$ \eqref{eqSigmaReal} and the following calculation,
\begin{align*}
iC^{[*]} J C = &
iC^{[*]} \left(I-I + J\right) C
\\
= &
iC^{[*]}C - 2iC^{[*]} Q C
\\
= &
(A-A^{[*]}) - iC^{[*]} Q C + (iC Q C)^{[*]}
\\
= &
(A - iC^{[*]} Q C) - (A - (iC Q C))^{[*]}
\\
= &
A^{\prime} - A^{\prime[*]}.
\end{align*}

\begin{Step}
We show $\mathcal P_K(S)$ does not contain nonzero constants.
\end{Step}
Let us consider the sequence
$ w_n \, k_{w_n}$ of $ \mathcal P (S)$
where $w_n = i \, n$ goes to infinity.
\smallskip

We show that the sequence has a limit in $\mathcal P (S)$ by applying Theorem \ref{IKL24}.
First, we note that, for any $\omega \in \Omega$ the following identity holds
\begin{align*}
\innerProductIndefiniteReg
{\frac{-J+S(\alpha)^{*} J S(z)}{-i(z - \overline{\alpha})}}
{\lim_{w \rightarrow \infty} w K_w(z)}
& =
\lim_{w \rightarrow \infty}
w
\frac{-J+S(\alpha)^{*} J S(w)}{-i(w - \overline{\alpha})}
\\
& =
i
\left(
-J + S(\alpha)^{*} J
\right)
\end{align*}
and thus the limit does exist.
Furthermore, a direct calculation yields
\begin{align}
\label{SecCalc}
&
\innerProductIndefiniteReg
{w_n K( \cdot , w_n) - w_m K( \cdot , w_m) }
{w_n K( \cdot , w_n) - w_m K( \cdot , w_m)}
= \\
\nonumber
&
|w_n|^2 K( w_n,w_n) + |w_m|^2 K( w_m,w_m)
-
\\
\nonumber
&
w_n K(w_n,w_m) \overline{w_m}
- w_m K(w_m,w_n) \overline{w_n}.
\end{align}
Since $S$ is analytic in a neighborhood of infinity,
we have the following expansion
\begin{equation}
\label{expanTheta}
S(w) = I + \frac{a}{w} + \frac{b}{w^2} + O(w^{-3}),
\end{equation}
where $a$ and $b$ are complex numbers.
Thus, using \eqref{expanTheta}, one rewrites the kernel as
\begin{align}
\label{expanKernel}
K(z,w) = & \frac{- J + S(w)^{*} J S(z)}{-i(z- \overline{w})}
\\
\nonumber
= &
\frac{1}{-i(z- \overline{w})}
\left\{
-J +
\left(I + \frac{a}{w} + \frac{b}{w^2}\right)^{*}
J
\left(I + \frac{a}{z} + \frac{b}{z^2}\right)
\right\}
\\
\nonumber
= &
\frac{1}{-i(z- \overline{w})}
\left(
J \frac{a}{z} +
\frac{a^*}{\overline  {w}} J
+
\cdots
\right)
.
\end{align}
Substituting \eqref{expanKernel} into \eqref{SecCalc}, leads to
\begin{align}
\nonumber
&
\innerProductIndefiniteReg
{w_n K( \cdot , w_n) - w_m K( \cdot , w_m) }
{w_n K( \cdot , w_n) - w_m K( \cdot , w_m)}
\\
\nonumber
=
&
\frac{J a \overline{w_n} + a^* J w_n}{-i(w_n - \overline{w_n})} +
\frac{J a \overline{w_m} + a^* J w_m}{-i(w_m - \overline{w_m})} -
\\
\nonumber
&
\frac{J a \overline{w_n} + a^* J w_m}{-i(w_m - \overline{w_n})} -
\frac{J a \overline{w_m} + a^* J w_n}{-i(w_n - \overline{w_m})}
+
O(w_n^{-1},w_m^{-1})
\\
\nonumber
=
&
i(J a + a^* J)
\frac
{|w_m-w_n|^2 \imag (w_m w_n)}
{-  2 \imag(w_m) \imag(w_n)|w_n - \overline{w_m}|^2}
+
O(w_n^{-1},w_m^{-1})
.
\end{align}
Since  $ w_n = i \, n$, the expression above tends to zero.
We then apply Theorem \ref{IKL24} and conclude the limit belongs to
$\mathcal P (S)$ and equal to $J- S (z)  J$.
Let $f(z)=c$ be a constant function in $\mathcal P (S)$. 
By the reproducing kernel property, we have that
\[
\innerProductReg{c}{i n K(\cdot, i \, n)\xi} =
-i \, n \overline{\xi} c.
\]
Taking $n \rightarrow \infty$, 
we have that the left hand side has a finite limit and so has the right hand side
and so $c$ equal zero.
It follows that $\mathcal P_K(S)$ does not contain nonzero constants.
\end{pf}

%%%%%%%%%%%%%%%%%%%%%%%%%%%%%%%%%%%%%
%									%
% Operator colligation over PS		%
%									%
%%%%%%%%%%%%%%%%%%%%%%%%%%%%%%%%%%%%%
\subsection{Single operator colligation over Pontryagin space}
We start with the definition of a single operator colligation over Pontryagin space.
\begin{definition}
A single operator colligation over Pontryagin spaces is given by the collection
\[
\mathcal C =
\left( A \, ; \, \mathcal P , \, \Phi \, , E \, ; \, \sigma \, \right),
\]
where $A$ is a linear non-selfadjoint operator defined on the Pontryagin space $\mathcal P$, such that
\begin{equation}
\label{collCond}
\frac{1}{i} \left( A - A^{[*]} \right) = \Phi^{[*]} \sigma \Phi.
\end{equation}
$E$ is a Pontryagin space, $\Phi$ is a linear mapping from $\mathcal P $ to $E$ and
$\sigma$ is a selfadjoint operator on $E$.
\end{definition}
Equation \eqref{collCond} is referred as the colligation condition of the operator $A$.
We use the notation $\sigma$ rather than $J$ since $\sigma$ is not necessarily a signature matrix.

The characteristic function is the single variable operator-valued function on $E$ defined by
\begin{equation}
\label{collSingleOp1Vars}
S(z) = I - i \Phi \left( A - z I \right) ^{-1} \Phi ^ {[*]} \sigma.
\end{equation}
The principal subspace is given by
\[
\widehat{\mathcal P}
=
\bigvee_{n=0}^{\infty}{A^n\Phi^{[*]}(E)}
=
\bigvee_{n=0}^{\infty}{A^{[*]n}\Phi^{[*]}(E)}.
\]
where $\bigvee$ denotes the closure of the union of the spaces $A^n\Phi^{[*]}(E)$ in $\mathcal P$.
It is worth noticing that for a neighborhood of infinity, say $\Omega$, the principal subspace is also given by
\[
\widehat{\mathcal P}
=
\bigvee_{z\in \Omega}^{}{(zI-A)^{-1}\Phi^{[*]}(E)}.
\]
It makes sense (see Remark \ref{degenerateColl} below) to introduce the following definitions.
\begin{definition}
A single-operator colligation over Pontryagin space $\mathcal P$ is called
\begin{enumerate}
\item
{\rm non--degenerate} if $\widehat{\mathcal P} \subset \mathcal P$ is a non--degenerate subspace (and hence Pontryagin subspace).
\item
{\rm irreducible} if $\widehat{\mathcal P} = \mathcal P$ (and in particular $\mathcal P$ is a Pontryagin space).
\end{enumerate}
\end{definition}

The description of the class of characteristic functions over Pontryagin spaces is given below.
\begin{theorem}
\label{singleOpColl}
An $n \times n$ matrix function $S(z)$ is the characteristic function of an irreducible operator colligation over Pontryagin space
of negative index $\kappa$ with invertible selfadjoint matrix $\sigma$ if and only if
\begin{enumerate}
\item
$S(z)$ is holomorphic in a neighborhood of $\infty$ and $S(\infty) = I$.
\item
$S(z)$ is meromorphic on the complement of the real axis, and for all $z,w$ in the region of analytically
\begin{equation}
\label{singleOpCollKernel}
\frac{S(w)^{*} \sigma S(z)  - \sigma}{-i (z - \overline{w})}
\end{equation}
has exactly $\kappa$ negative squares.
\end{enumerate}
\end{theorem}
The counterpart theorem in the Hilbert case can be found in \cite[Theorem 5.1]{brod}.
\begin{Rk}
If we endow $E$ with the Hermitian form
\begin{equation}
\label{sigmaIdefiniteE}
\innerProductReg{x}{\sigma y} = \innerProductIndefiniteTri{x}{y}{\sigma},
\end{equation}
then \eqref{singleOpCollKernel} becomes
\[
\frac{S(w)^{[*]}S(z)-I}{-i(z-\overline{w})},
\]
where $[*]$ is the adjoint with respect to \eqref{sigmaIdefiniteE}.
By a generalization of \eqref{adjoint2adjoint} to Hermitian matrix case, i.e.
\[
A^{[*]} = \sigma^{-1} A^* \sigma,
\]
it follows that
\begin{align*}
\innerProductIndefiniteTri{S(w)^{[*]}S(z)x}{y}{\sigma}
= &
\innerProductIndefiniteTri{S(z)x}{S(w) y}{\sigma}
\\ = &
y^* \sigma S(w)^{[*]} S(z) x
\\ = &
y^*  S(w)^{*} \sigma S(z) x.
\end{align*}
\end{Rk}

\begin{pf}[of Theorem \ref{singleOpColl}]
Let $\mathcal C$ be an irreducible single operator colligation with characteristic function $S(z)$ of the form \eqref{collSingleOp1Vars}.
Clearly, it is holomorphic in the neighborhood of infinity, $S(\infty)=I$ and $S(\overline{z})^* \sigma S(z) = \sigma$.
Furthermore,
using the colligation condition and applying a well-known classical calculation,
\begin{align*}
\nonumber
K_w(z)
= &
\frac
{- \sigma + (I  - i  \sigma^* \Phi (wI - A)^{-[*]} \Phi^{[*]} ) \sigma (I + i \Phi (zI - A)^{-1} \Phi^{[*]} \sigma)}
{-i(z-\overline{w})}
\\
\nonumber
= &
i
\frac
{
  \sigma \Phi (zI - A)^{-1} \Phi^{[*]} \sigma
    - \sigma^* \Phi (wI - A)^{-[*]} \Phi^{[*]} \sigma
}
{-i(z-\overline{w})}
+
\\ & \nonumber + 
\frac
{
\sigma^* \Phi (wI - A)^{-[*]} \Phi^{[*]}  \sigma  \Phi (zI - A)^{-1} \Phi^{[*]} \sigma
}
{-i(z-\overline{w})}
\\
= &
i
\sigma^* \Phi   (wI - A)^{-[*]} 
\frac
{
 A^{[*]} - A - i \Phi^{[*]} \sigma \Phi + (z - \overline{w})
}
{-i(z-\overline{w})}
    (zI - A)^{-1}\Phi^{[*]} \sigma
,
\end{align*}
one may conclude that
\[
S(w)^{*}\sigma S(z) - \sigma
=
i
\sigma^* \Phi   (wI - A)^{-[*]} 
(zI - A)^{-1}\Phi^{[*]} \sigma
.
\]
Thus, using Theorem \ref{pontKernelDecomA}, \eqref{singleOpCollKernel} has exactly $\kappa$ negative squares.
\smallskip

Assuming the converse, let $S(z)$ be holomorphic at a neighborhood of infinity such that \eqref{singleOpCollKernel} has $\kappa$ negative squares.
Using Theorem \ref{th1}, there exists a colligation
\[
\mathcal C_{S} = \left(A \, ; \, \mathcal P_K (  S ) \, , \, \Phi \, , \, \mathbb C ^n \, ; \,  \sigma  \,   \right),
\]
such that $S(z)$ is its characteristic function.
Here, $\mathcal P _K(  S )$ is the reproducing kernel Pontryagin space with reproducing kernel \eqref{singleOpCollKernel},
$A$ is a bounded linear operator on $\mathcal P_K (  S )$,
$E$ is the Pontryagin space $\mathbb C ^n$ equipped with the indefinite inner product $\left[u,v\right]_{\sigma} = v^{[*]} \sigma u$.
Note that the colligation condition is just \eqref{eqCollForS}.
\smallskip

Finally, since the $ \lim_{z \rightarrow \infty} {S (z)} = I$, one may conclude that $\mathcal P_K (S)$ does not contain nonzero constant functions.
\end{pf}

%%%%%%%%%%%%%%%%%%%%%%%%%%%%%%%%%%%%%
%									%
%	Pontryagin vessel theory		%
%									%
%%%%%%%%%%%%%%%%%%%%%%%%%%%%%%%%%%%%%
\section{Vessel theory over Pontryagin spaces}
\setcounter{equation}{0}

Vessel theory over Hilbert spaces is well developed,
see for instance \cite{MR595748,livsic2,livsic1,MR1634421} and the book by Liv\v{s}ic et al. \cite{KLMV}.
In this section we set the counterpart theory for the non-positive setting.

%%%%%%%%%%%%%%%%%%%%%%%%%%%%%%%%%%%%%
%									%
%	Colligations					%
%									%
%%%%%%%%%%%%%%%%%%%%%%%%%%%%%%%%%%%%%
\subsection{The vessel definition and some basic properties}
\begin{definition}
A {\it commutative two-operator colligation over Pontryagin space $\mathcal P$} is the collection
\[
\mathcal C
\defEq
\left(
 \, A_1 \, , \, A_2 \,  ; \,  \mathcal P \,  , \,  \Phi \,  , \,  E \,  ; \,  \sigma_1 \, ,  \, \sigma_2 \,
\right)
,
\]
where the non-selfadjoint operators $A_1$ and $A_2$ commute ($A_1 A_2 = A_2 A_1$) and satisfy
\begin{equation}
\label{vesselCollCond}
\frac{1}{i}\left(A_k - A_k^{[*]} \right) =\Phi^{[*]} \sigma_k \Phi, \qquad k=1,2.
\end{equation}
In this setting,
$E$, the outer space, is a finite dimensional Pontryagin space,
$\Phi$ is a linear mapping from $\mathcal P$ to $E$
and $\sigma_1$ and $\sigma_2$ are selfadjoint operators in $E$.
\end{definition}
Equations \eqref{vesselCollCond} are referred as the colligation conditions associated to $A_1$ and $A_2$.
It is well-known, see for instance \cite{lifshits1},
that a two-operator colligation does not contain enough information regarding the interplay between the two operators.
This observation naturally leads to the definition of vessels over Pontryagin spaces.
\begin{definition}
The collection
\[
\mathcal V =
\left(
 \, A_1 \, , \,  A_2  \, ;  \, \mathcal P \, , \,  \Phi \, ,  \, E  \, ;  \, \sigma_1 \, ,  \, \sigma_2 \,  , \,  \gamma  \, , \, \widetilde{\gamma}  \,
\right)
\]
is called a commutative two-operator vessel over Pontryagin space if
\allowbreak
$\left(A_1,A_2;\mathcal{P},\Phi,E;\sigma_1,\sigma_2 \right)$
is a commutative two-operator colligation,
$\gamma$ and $\widetilde{\gamma}$ are selfadjoint operators on $E$
and furthermore the following conditions hold:
\begin{align}
\label{vesselPSinput}
\gamma \Phi & = \sigma_1 \Phi A_2^{[*]} - \sigma_2 \Phi A_1^{[*]}
\\
\label{vesselPSoutput}
\widetilde{\gamma} \Phi & = \sigma_1 \Phi A_2 - \sigma_2 \Phi A_1
\\
\label{vesselPSlinkage}
\widetilde{\gamma} - \gamma
& =
i
\left(
\sigma_1 \Phi \Phi^{[*]} \sigma_2 - \sigma_2 \Phi \Phi^{[*]} \sigma_1
\right)
.
\end{align}
\end{definition}
Equations \eqref{vesselPSinput}, \eqref{vesselPSoutput} and \eqref{vesselPSlinkage} are known as the
input, output and linkage conditions of the vessel, respectively.
It is worth mentioning that condition \eqref{vesselPSinput} follows
from \eqref{vesselCollCond}, \eqref{vesselPSoutput} and \eqref{vesselPSlinkage}
(the same observation holds for \eqref{vesselPSoutput}).
\subsection{The principal subspace}
\begin{theorem}
Let $\mathcal C$ be a two-operator commutative colligation over Pontryagin space. Then the following relation holds
\[
\widehat{\mathcal P}
\defEq
\bigvee_{n=0}^{\infty}{A_1^{k_1}A_2^{k_2}\Phi^{[*]}(E)}
=
\widehat{\mathcal P}^{[*]}
\defEq
\bigvee_{n=0}^{\infty}{A_1^{[*] k_1}A_2^{[*] k_2}\Phi^{[*]}(E)}.
\]
\end{theorem}
\begin{pf}
Let $f \in \widehat{\mathcal P}$ and let us assume that
\[
A_1^{k_1} A_2^{k_2} \Phi^{[*]}(E) \subset  \widehat{\mathcal P}^{[*]}.
\]
We note that
\[
A_1^{k_1+1} A_2^{k_2} \Phi^{[*]}(E)
=
(A_1 - A_1^{[*]})A_1^{k_1} A_2^{k_2} \Phi^{[*]}(E) +
A_1^{[*]}A_1^{k_1} A_2^{k_2} \Phi^{[*]}(E),
\]
by the colligation condition
\[
A_k f  - A_k^{[*]}f \in \widehat{\mathcal P}^{[*]}.
\]
Hence (and by the fact that $\widehat{\mathcal P}$ is invariant under $A_1^{[*]}$)
\[
A_1^{k_1+1} A_2^{k_2} \Phi^{[*]}(E)
\subset \widehat{\mathcal P}^{[*]}.
\]
By induction on $k_1$ and $k_2$ we may conclude that $\widehat{\mathcal P} \subseteq \widehat{\mathcal P}^{[*]}$.
Similarly, $\widehat{\mathcal P}^{[*]} \subseteq \widehat{\mathcal P}$ and hence $\widehat{\mathcal P}^{[*]} = \widehat{\mathcal P} $.
\end{pf}
\begin{Rk}
\label{degenerateColl}
The subspace $\widehat{\mathcal P}$ is called the {\it principal subspace} of two-operator colligation $\mathcal C$.
In opposition to the Hilbert space case $\widehat{\mathcal P}$ may be degenerate.
In view of Theorems \ref{bognar22} and \ref{bognar23}, we are interested in the case where the principal subspace $\widehat{\mathcal P}$ is also a Pontryagin space.
This observation leads us to focus in the case where the principal subspace is non--degenerate. 
\end{Rk}

\begin{definition}
A vessel over Pontryagin space $\mathcal P$ is called
\begin{enumerate}
\item
{\rm non--degenerate} if $\widehat{\mathcal P} \subseteq \mathcal P$ is non--degenerate subspace (and hence Pontryagin subspace).
\item
{\rm irreducible} if $\widehat{\mathcal P} = \mathcal P$ (and in particular $\mathcal P$ is a Pontryagin space).
\end{enumerate}
\end{definition}
Clearly, the principal subspace of a vessel is invariant under the operators $A_1$ and $A_2$
(and using the colligation condition is also invariant under their ajoints $A_1^{[*]}$ and $A_2^{[*]}$).
Thus, it is also true that the orthogonal complement of $\widehat{\mathcal P}$ is invariant under $A_1$, $A_2$, $A_1^{[*]}$ and $A_2^{[*]}$.
\smallskip

Furthermore, the principal subspace a vessel contains the whole non-Hermitian structure of these operators.
In other words, the restrictions of $A_1$ and $A_2$ to $\mathcal P \ominus \widehat{\mathcal P}$ are selfadjoint operators. 
It is true, since for $v \in \mathcal P \ominus \widehat{\mathcal P}$, the following identity holds
\[
0
=
\innerProductIndefiniteTri{v}{\Phi^{[*]}E}{\mathcal P}
=
\innerProductIndefiniteTri{\Phi v}{E}{E}.
\]
Hence $\Phi v = 0$, and by the colligation condition, we have the following
\[
(A_k - A_k^{[*]})v = i \Phi^{[*]} \sigma_k \Phi v = i \Phi^{[*]} \sigma_k 0 = 0.
\]
Thus, $A_k \restrict{\mathcal P \ominus \widehat{\mathcal P}}$ are selfadjoint operators.
\smallskip

Two vessels are said to be unitary equivalent if they have common $E,\sigma_1,\sigma_2$ and
there exists a unitary operator from $\widehat{\mathcal P _1}$ to $\widehat{\mathcal P _2}$ such that
\[
A_k^2 \restrict{\widehat{\mathcal P _2}} =  U A_k^1 \restrict{ \widehat{\mathcal P _2}} U^{-1},
\qquad
\Phi^2\restrict{\widehat{\mathcal P _2}} = \Phi^1 \restrict{ \widehat{\mathcal P _1}} U,
\]
where $k=,1,2$.
In particular, their principal subspaces share the same negative index
and have the same complete characteristic function.
%%%%%%%%%%%%%%%%%%%%%%%%%%%%%%%%%%%%%
%									%
%	C-H theorem						%
%									%
%%%%%%%%%%%%%%%%%%%%%%%%%%%%%%%%%%%%%
\subsection{The complete characteristic function}
The {\it complete characteristic function} (CCF) of a vessel is defined by
\[
W(\xi_1,\xi_2,z) = I - i \Phi( \xi_1 A_1 + \xi_2 A_2 - zI)^{-1}\Phi^{[*]}(\xi_1 \sigma_1 + \xi_2 \sigma_2),
\]
where $\xi_1,\xi_2,z \in \mathbb C$.
It is a homogenous function of degree zero and analytic function for all
$z \in \mathbb C$ outside the spectrum of $\xi_1 A_1 + \xi_2 A_2$.
The CCF, as a function of several variables, does not admit a good factorization theory.
However, it has some striking properties.
\begin{Rk}[The $\alpha$-transformation]
\label{alphaTrans}
The $\alpha$-transformation of a vessel $\mathcal V$ (see also \cite{livsic1} or \cite{KLMV})
is the following linear transformation
\begin{align*}
\ddot A_1 & = \alpha_{11} A_1 + \alpha_{12}A_2, \quad
\ddot{\sigma}_1  = \alpha_{11} \sigma_1 + \alpha_{12} \sigma_2, \quad
%z^{\prime}_1 & = \alpha_{11} z_1 + \alpha_{12} z_2, \quad
\\
\ddot A_2 & = \alpha_{21} A_1 + \alpha_{22}A_2, \quad
\ddot{\sigma}_2  = \alpha_{11} \sigma_1 + \alpha_{22} \sigma_2, \quad
% z^{\prime}_2 & = \alpha_{21} z_1 + \alpha_{22} z_2,
\end{align*}
    where $\alpha \in \mathbb C ^{2 \times 2}$ with $\det{\alpha}=1$.
Under the $\alpha$ transformation, the collection
\[
\ddot{ \mathcal V}=
\left(
\ddot A_1 \, , \,  \ddot A_2  \, ; \,   \mathcal P \, , \,  \Phi \, , \,  E \,  ; \,  \ddot \sigma_1 \, , \,  \ddot \sigma_2  \, , \,  \gamma \,  , \, \widetilde{\gamma}
\right)
\]
is a Pontryagin space vessel with the complete characteristic function
\[
\ddot{W}(\xi_1,\xi_2,z) = W(\alpha_{11} \xi_1 + \alpha_{12} \xi_2,\alpha_{21} \xi_1 + \alpha_{22} \xi_2,z).
\]
\end{Rk}
It is useful to define to following subspaces of $E$:
\begin{align*}
E(\lambda_1,\lambda_2) & = \ker (\lambda_1 \sigma_2 - \lambda_2 \sigma_1 +\gamma)\\
\widetilde{E}(\lambda_1,\lambda_2) & = \ker (\lambda_1 \sigma_2 - \lambda_2 \sigma_1 +\widetilde{\gamma}).
\end{align*}
\begin{lem}
\label{jcfLem}
The complete characteristic function
$W(\xi_1,\xi_2,\xi_1 \lambda_1 + \xi_2 \lambda_2)$
defines a mapping from
$E(\lambda_1,\lambda_2)$ to $\widetilde{E}(\lambda_1,\lambda_2)$.
Furthermore, its restriction to $E(\lambda_1,\lambda_2)$ is independent of the choice of $\xi_1$ and $\xi_2$.
\end{lem}
\begin{pf}
Let $e$ be an element in $E(\lambda_1,\lambda_2)$. Then, by definition, we have $(\lambda_1\sigma_1 - \lambda_2 \sigma_2 +\gamma)e = 0$
and hence $\Phi^{[*]}(\lambda_1\sigma_1 - \lambda_2 \sigma_2 +\gamma)e = 0$.
Now using \eqref{vesselPSinput} we may conclude that
\begin{equation}
\label{ccfProp1}
(A_1-\lambda_1 I)\Phi^{[*]} \sigma_2 e = (A_2-\lambda_2 I)\Phi^{[*]} \sigma_1 e.
\end{equation}
Using the $\alpha$-transformation (see Remark \ref{alphaTrans}), \eqref{ccfProp1} becomes
\[
(A^{\prime}_1-\lambda^{\prime}_1 I)\Phi^{[*]} \sigma^{\prime}_2 e = (A^{\prime}_2-\lambda^{\prime}_2 I)\Phi^{[*]} \sigma^{\prime}_1 e.
\]
Assuming that $\xi_1 \lambda_1 + \xi_2 \lambda_2$ does not belong to the spectrum of $\xi_1 A_1 + \xi_2 A_2$,
then $\eta_1 \lambda_1 + \eta_2 \lambda_2$ does not belong to the spectrum of $\eta_1 A_1 + \eta_2 A_2$
and using the commutativity of $A_1$ and $A_2$. Thus
\[
(\xi A - \xi \lambda)^{-1} \Phi^{[*]} (\xi_1 \sigma_1 + \xi_2 \sigma_2) e.
\]
is independent of the choice of $\xi_1$ and $\xi_2$ and so does $W(\xi_1,\xi_2,\xi_1 y_1 + \xi_2 y_2)$.
To complete the proof we note that \eqref{eq33} indicates that
$$W(\xi_1,\xi_2,\xi_1 \lambda_1 + \xi_2 \lambda_2)$$ maps $E(\lambda_1,\lambda_2)$
to $\widetilde{E}(\lambda_1,\lambda_2)$.
\end{pf}
\subsection{The discriminant variety and the joint characteristic function}
We continue and show that the input and output discriminant polynomials coincide 
(the input and the output terms are coined due to system theory interpretation).
\begin{theorem}
Let $\mathcal V$ be a commutative Pontryagin vessel.
Then the input and the output discriminant polynomials coincide
    \begin{equation}
        \label{eqDetRep}
p(z_1,z_2)=
\det \left( z_1 \sigma_2 - z_2 \sigma_1 + \gamma \right) =
\det \left( z_1 \sigma_2 - z_2 \sigma_1 + \widetilde{\gamma} \right).
    \end{equation}
\end{theorem}
\begin{pf}
Using \eqref{vesselPSinput} and \eqref{vesselPSoutput}, we have
\begin{align*}
(A_2 - z_2 I) \Phi ^{[*]} \sigma_1
-
(A_1 - z_1 I) \Phi ^{[*]} \sigma_2
= &
\Phi ^{[*]}
\left(
z_1 \sigma_2 - z_2 \sigma_1 + \gamma
\right)
\\
\sigma_1 \Phi ^{[*]} (A_2 - z_2 I)
-
\sigma_2 \Phi ^{[*]} (A_1 - z_1 I)
= &
\left(
z_1 \sigma_2 - z_2 \sigma_1 + \widetilde{\gamma}
\right)
\Phi.
\end{align*}
Multiplying the first equality by $\sigma_1 \Phi(A_1 - z_1 I)^{-1}$ on the left and the second by $\sigma_2 \Phi(A_2 - z_2 I)^{-1}$ on the right we have
\begin{align*}
\sigma_2 \Phi \Phi^{[*]} \sigma_1 - \sigma_1 \Phi \Phi^{[*]} \sigma_2
=
&
\sigma_1 \Phi ^{[*]} (A_1 - z_1 I)^{-1}
\Phi ^{[*]}
\left(
z_1 \sigma_2 - z_2 \sigma_1 + \gamma
\right)
\\
&
-
\left(
z_1 \sigma_2 - z_2 \sigma_1 + \widetilde{\gamma}
\right)
\Phi (A_1 - z_1 I)^{-1}  \Phi ^{[*]}  \sigma_1.
\end{align*}
Now using \eqref{vesselPSlinkage} we can conclude that
\begin{equation}
\label{eq33}
\left(
z_1 \sigma_2 - z_2 \sigma_1 + \widetilde{\gamma}
\right)
W(1 , 0 ,  z_1 )
=
\widetilde{W}(1 , 0 ,  z_1 )
\left(
z_1 \sigma_2 - z_2 \sigma_1 + \gamma
\right)
\end{equation}
where
\[
    \widetilde{W}(\xi_1 , \xi_2 ,  z ) = I - i(\xi_1 \sigma_1 + \xi_2 \sigma_2)\Phi(\xi_1 A_1 + \xi_2 A_2 -zI)^{-1} \Phi^{[*]}.
\]
We note that for $R$ large enough, $\det {W(1 , 0 ,  z_1) } \neq 0 $ where $|z_1|>R$ and by using $\det (I + M_1 M_2) = \det (I + M_2 M_1)$, one may obtain
\[
\det W(1,0,z_1) = \det \widetilde{W} (1,0,z_1).
\]
Thus by \eqref{eq33}
\[
\det \left( z_1 \sigma_2 - z_2 \sigma_1 + \gamma \right) =
\det \left( z_1 \sigma_2 - z_2 \sigma_1 + \widetilde{\gamma} \right).
\]
\end{pf}

The polynomial $p(\lambda_1,\lambda_2)$ is called the {\it discriminant polynomial} of the vessel $\mathcal V$.
Hence a commutative vessel defines an affine discriminant curve,
\[
C_0 \defEq \big\{ (\lambda_1,\lambda_2)\in \mathbb R ^2 \quad | \quad \det \left( \lambda_1 \sigma_2 - \lambda_2 \sigma_1 + \gamma \right) = 0 \big\},
\]
with two  associated selfadjoint determinantal representations (see \eqref{eqDetRep}).
\smallskip

Using Lemma \ref{jcfLem}, we now define a function of a single variable on the discriminant curve.
\begin{definition}
The {\it joint characteristic function} of a vessel over Pontryagin space $\mathcal P$ is defined by
\[
S(\lambda_1,\lambda_2) \defEq
\restr{W(\xi_1,\xi_2,\xi_1 \lambda_1 + \xi_2 \lambda_2)}{E(\lambda_1,\lambda_2)} : E(\lambda_1,\lambda_2) \rightarrow \widetilde{E}(\lambda_1,\lambda_2),
\]
where $(\lambda_1,\lambda_2)$ belongs to $C_0$ and $\xi_1 \lambda_1 + \xi_2 \lambda_2$ does not belong to the spectrum of $\xi_1 A_1 + \xi_2 A_2$.
\end{definition}

The following theorem is purely algebraic and has the same proof as in the Hilbert space case
as long as one assumes that the principal subspace is non-degenerate.
\begin{theorem}[Generalized Cayley-Hamilton theorem]
Let $\mathcal V$ be a non--degenerate vessel over Pontryagin space.
Then following identity holds
\[
p(A_1,A_2)
=
0
\]
on its principal subspace.
\end{theorem}
Finally, we conclude this section with an important result used later in this sequel
(see Section \ref{sec2dSys} and Section \ref{secRealization}).
\begin{lem}
\label{lemmaPrincSubSpace}
Let $\mathcal V$ be a commutative two-operator vessel over Pontryagin space.
\begin{enumerate}
\item
Let $\Omega$ be a neighborhood of $0 \in \mathbb R ^2$, then
\begin{align}
\mathcal P
\nonumber
= &
\bigvee_{(t_1,t_2)\in \Omega}^{\infty}{\exp{(it_1A_1+it_2A_2)}\Phi^{[*]}(E)} \\
\label{princSubSpace}
= &
    \bigvee_{(t_1,t_2)\in \Omega}^{\infty}{\exp{(it_1A_1^{[*]}+it_2A_2^{[*]})}\Phi^{[*]}(E)}.
\end{align}
\item
\label{princSubSpace2}
Furthermore, assume $\mathcal V$ is non--degenerate.
If
\[
\Phi \exp{(it_1A_1+it_2A_2)}p=0,
\qquad
\forall (t_1,t_2)\in \Omega
\]
where $p\in \widehat{\mathcal P}$, then $p = 0$.
\end{enumerate}
\end{lem}
\begin{pf}
For the first part, see \cite[Proposition 10.4.1 (1)]{KLMV}.
The second part is similar to \cite[Proposition 10.4.1 (2)]{KLMV}, but note that the assumption that $\mathcal V$ is non-degenerate is crucial. 
Since $\Phi \exp{(it_1A_1+it_2A_2)}p = 0$ for every $(t_1,t_2)\in W$, by analyticity it is also true for all $(t_1,t_2)\in \mathbb R^2$.
We set $$G = \bigvee_{(t_1,t_2)\in \mathbb R^2}\exp{(it_1A_1+it_2A_2)}p.$$
Then, for any $e \in E$ and $(t_1,t_2)\in \mathbb R^2$, we have
\[
\innerProductIndefiniteTri{\exp{(it_1A_1+it_2A_2)}h}{\Phi^{[*]}e}{\mathcal P} =
\innerProductIndefiniteTri{\Phi \exp{(it_1A_1+it_2A_2)}h}{e}{E} = 0.
\]
Thus $G \perp \Phi^{[*]} E$, and so $\Phi^{[*]} E \subset G ^\perp$ and
$\exp{(it_1A_1+it_2A_2)}\Phi^{[*]} E \subset G ^\perp$ and then
by \eqref{princSubSpace} we have  $\widehat{\mathcal P} \subset G ^\perp$.
Thus $p = 0$ since $p \in \widehat{\mathcal P}$ and $\widehat{\mathcal P}$ is non--degenerate.
\end{pf}

%%%%%%%%%%%%%%%%%%%%%%%%%%%%%%%%%%%%%
%									%
% 				$2D$ systems		%
%									%
%%%%%%%%%%%%%%%%%%%%%%%%%%%%%%%%%%%%%
\section{Overdetermined \texorpdfstring{$2D$}{ s } systems over indefinite inner spaces}
\label{sec2dSys}
\setcounter{equation}{0}
A main motivation to study vessel theory over indefinite inner product spaces
is the study of overdetermined $2D$ systems in the non-positive setting.
A general $2D$ linear time--invariant system is given by the set of equations:
\begin{align}
\label{2DsysA}
i \frac{\partial }{\partial t_1}f(t_1,t_2) + A_1 f(t_1,t_2) & = B_1 u(t_1,t_2),
\\
\label{2DsysB}
i \frac{\partial}{\partial t_2}f(t_1,t_2) + A_2 f(t_1,t_2) & = B_2 u(t_1,t_2),
\\
\label{2DsysC}
v(t_1,t_2) & = u(t_1,t_2) + C f(t_1,t_2).
\end{align}
The functions $u(t_1,t_2)$, $f(t_1,t_2)$ and $v(t_1,t_2)$ are functions of two variables $t_1$ and $t_2$ denoted by the input, system and output function of the system, respectively.
The operators $A_1$ and $A_2$ are bounded in a Pontryagin space $\mathcal P$ with $\kappa$ negative index.
$E$ is a finite dimensional Pontryagin space, $C:\mathcal P \rightarrow E$ and $B: E \rightarrow \mathcal P$.
\smallskip

We assume that the partial derivatives commute, that is
$\frac{\partial^2}{\partial t_1 \partial t_2} =\frac{\partial^2}{\partial t_2 \partial t_1}$.
As a consequence, \eqref{2DsysA} and \eqref{2DsysB} become
\begin{equation}
\label{2DsysD}
\left[
\left(
B_1 \frac{\partial}{\partial t_2} - B_2 \frac{\partial}{\partial t_1}
\right)
+i
\left(
A_1 B_2 - A_2 B_1
\right)
\right]
u
=
i
\left(
A_1 A_2 - A_2 A_1
\right)
f
.
\end{equation}

Furthermore, we assume that the system preserves the scattering energy balance (see \cite{livvsic1987commuting})
\[
\frac{d}{dt_k} \innerProductIndefiniteTri{f}{f}{\mathcal P}
=
\innerProductIndefiniteTri{\sigma_k u}{u}{E} - \innerProductIndefiniteTri{\sigma_k y}{y}{E},
\qquad
k=1,2,
\]
that is, the change of the state energy with respect to each variable is equal to the difference between the output energy and the input energy.
Then, comparing (\ref{2DsysA}-\ref{2DsysC}) with the scattering energy balance equations, we have
\begin{align*}
\innerProductIndefiniteReg{(A_k - A_k^{[*]})f}{f}
+
&
\innerProductIndefiniteReg{f}{B_1u}
-
\innerProductIndefiniteReg{B_1u}{f}
=
\\
\nonumber
&
\innerProductIndefiniteReg{iC^{[*]}\sigma_k C f}{f}
+
\innerProductIndefiniteReg{i \sigma_k C f}{u}
+
\innerProductIndefiniteReg{u}{-\sigma_k C f}.
\end{align*}
Assuming consistency for zero input, $u(t_1,t_2)=0$, and using \eqref{2DsysD}, we obtain that the operators $A_1$ and $A_2$ commute (see also \cite{MR2043236}).
By the same consistency assumption, we have two additional immediate results. The first (the colligation condition):
\[
C^{[*]} \sigma_ k C = \frac{1}{i} ( A_k - A_k^{[*]}), \qquad k=1,2
\]
and the second is $B_k =  - i C^{[*]} \sigma_k$ for $k=1,2$.
Hence under energy balance condition and the commutativity of the mixed derivatives, we consider the following overdetermined $2D$ system.
\begin{definition}
A $2D$ overdetermined linear time--invariant system over Pontryagin space is the collection
\[
\mathcal C =
\left( A_1 \, , \, A_2 ; \, \mathcal P , \, \Phi \, , E \, ; \, \sigma_1 \, , \, \sigma_2 \, \right),
\]
satisfying the following equations:
\begin{align*}
i \frac{\partial}{\partial t_1}(t_1,t_2) + A_1 f(t_1,t_2) & = \Phi^{[*]} \sigma_1 u(t_1,t_2),
\\
i \frac{\partial}{\partial t_2}f(t_1,t_2) + A_2 f(t_1,t_2) & = \Phi^{[*]} \sigma_2 u(t_1,t_2),
\\
v(t_1,t_2) & = u(t_1,t_2) - i \Phi f(t_1,t_2).
\end{align*}
Here $A_1$ and $A_2$ are commuting bounded operators in the Pontryagin space $\mathcal P$ with negative index $\kappa$, such that $A_k - A_k^{[*]} = i \Phi^{[*]}\sigma_k \Phi$ for $k=1,2$.
The functions $u$, $f$ and $v$ are two-variable functions of $t_1$ and $t_2$ denoted by the input, system and output function of the system, respectively.
$E$ is a finite dimensional Pontryagin space, $\sigma_k$ are selfadjoint operators on $E$ and $\Phi:\mathcal P \rightarrow E$.
\end{definition}
Furthermore, the compatibility conditions imply that,
for some selfadjoint operators $\gamma$ and $\widetilde{\gamma}$,
the following conditions,
\begin{align}
\label{compCondInput}
\left( \sigma_2 \frac{\partial}{\partial t_1} - \sigma_1 \frac{\partial}{\partial t_2} + i \gamma \right) u  & = 0
\\
\label{compCondOutput}
\left( \sigma_2 \frac{\partial}{\partial t_1} - \sigma_1 \frac{\partial}{\partial t_2} + i \widetilde{\gamma} \right) v  & = 0.
\end{align}
hold for any $u,v\in E$.
In fact, one can see that conditions \eqref{compCondInput} and \eqref{compCondOutput} are equivalent to conditions
\eqref{vesselPSinput}, \eqref{vesselPSoutput} and \eqref{vesselPSlinkage}.
Thus, commutative two-operator vessels over Pontryagin space are exactly the
overdetermined $2D$ systems over Pontryagin spaces equipped with the compatibility conditions \eqref{compCondInput} and \eqref{compCondOutput}.
This point is summarized in the following theorem.
\begin{Tm}
% \label{lemVesselEqSys2D}
Let $\mathcal V$ be a non--degenerate two-operator vessel over Pontryagin space (with negative index $\kappa$) of the form
\[
\mathcal V =
\left(
A_1 \, , \, A_2 \, ; \,  \mathcal P \, , \,  \Phi \, ,  \, E \,  ;
\,  \sigma_1 \, , \,  \sigma_2  \, , \,  \gamma \,  , \, \widetilde{\gamma}
\right)
.
\]
Then $\mathcal V$ is unitary equivalent on $\widehat{\mathcal P}$ to the irreducible two-operator vessel
\[
\widetilde{\mathcal V} =
\left(
-i \frac{\partial}{\partial t_1} \, , \,  -i \frac{\partial}{\partial t_2}  \, ;  \,
 \,  \ddot{\mathcal P} \, , \,  \Psi \,  , \,  E  \, ;  \, \sigma_1 \, , \,  \sigma_2 \,  , \,  \gamma \,  , \, \widetilde{\gamma}
\right)
,
\]
where $\ddot{\mathcal P}$ is Pontryagin space with ($\ddot{\kappa} \leq \kappa$) negative index containing solutions of the PDE
\[
\sigma_2 \frac{\partial}{\partial t_1}v -
\sigma_1 \frac{\partial}{\partial t_2}v +
i \widetilde{\gamma} v = 0.
\]
\end{Tm}
\begin{pf}
We note that since $\mathcal P$ is non--degenerate (and closed) then, by Theorem \ref{bognar23},
$\widehat{\mathcal P}$ is a Pontryagin space with
$\ddot{\kappa} (\leq \kappa)$ negative index.
\smallskip

Consider the $2D$ system corresponding to $\mathcal V$ and
assume the input $u(t_1,t_2)=0$ and the initial condition of the state is given by $f(0,0)=h$ where $h \in \widehat{\mathcal P}$.
Then the state is of the form $f_h(t_1,t_2) = \exp{(it_1A_1+it_2A_2)}h$ and the output is
$v_h(t_1,t_2) = -i \Phi \exp{(it_1A_1+it_2A_2)}h$.
\smallskip

Let $\ddot{\mathcal P} = \{v_p(t_1,t_2) \, | \, p \in \widehat{\mathcal P}\}$.
By Lemma \ref{lemmaPrincSubSpace} \eqref{princSubSpace2} the mapping
$U : \widehat{\mathcal P} \rightarrow \ddot{\mathcal P}$ given by $U p = v_p$ is injective.
Thus, $\ddot{\mathcal P}$ is Pontryagin space with $\ddot{\kappa}$ negative square with the indefinite inner product inherited from $\mathcal P$ and $U$ becomes an isometry.
Under this isometry $U$, the operators $A_1, A_2$ and $\Psi$ become
$-i \frac{\partial}{\partial t_1}, -i \frac{\partial}{\partial t_2}$
and $\Psi f(t_1,t_2) = i v_f(0,0)$, respectively.
\end{pf}

\smallskip

The relation between the input to the output signals of the $2D$ systems in the frequency domain is given below in term of the joint characteristic function of the system.
Let us assume that the input, output and state signals are waves with the same double period frequencies $(\lambda_1,\lambda_2) \in \mathbb C ^2$, that is
\begin{align*}
u(t_1,t_2) & = \widehat{u} e^{i t_1 \lambda_1 + i t_2 \lambda_2}, \qquad \widehat{u} \in E \\
v(t_1,t_2) & = \widehat{v} e^{i t_1 \lambda_1 + i t_2 \lambda_2}, \qquad \widehat{v} \in E \\
x(t_1,t_2) & = \widehat{x} e^{i t_1 \lambda_1 + i t_2 \lambda_2}, \qquad \widehat{x} \in \mathcal P.
\end{align*}
We note that condition \eqref{compCondInput} implies that $(\lambda_1,\lambda_2) \in C_0$ and $\widehat{u} \in E(\lambda)$
and similarly, condition \eqref{compCondOutput} implies that $(\lambda_1,\lambda_2) \in C_0$ and $\widehat{v} \in \widetilde{ E} (\lambda)$.
Then, the relation between the input and the output double periodic waves are given by
\[
\widehat{v} = S(\lambda) \widehat{u}, \qquad \lambda \in C_0.
\]
%%%%%%%%%%%%%%%%%%%%%%%%%%%%%%%%%%%%
%									%
%	Coupling $ decompositions		%
%									%
%%%%%%%%%%%%%%%%%%%%%%%%%%%%%%%%%%%%
\section{Pontryagin vessels \texorpdfstring{$2D$}{ s } systems coupling and decompositions}

The main purpose of this section is to study cascades and decompositions of
$2D$ systems in the non-positive case.
In the Hilbert space case, a closed subspace is
automatically ortho-complemented.
This fact is not true in the indefinite metric
and, as a consequence, we are forced to add the ortho-complemented hypothesis in the appropriate claims.

In order to decompose a system, we assume the inner space has a proper subspace,
which is a Pontryagin space (i.e. a closed and non--degenerate subspace of $\mathcal P$, see Theorems \ref{IKLthm5}-\ref{IKLthm6}) and
invariant under the two operators of the vessel.

Let $\mathcal V$ be a commutative two-operator vessel over a Pontryagin space such that there is a ortho--complemented closed subspace
$\mathcal P^2 \subseteq \mathcal P$ which is invariant under both operators $A_1$ and $A_2$.
Then we set $\mathcal P^1 = \mathcal P \ominus \mathcal P^2$ and cosidering the corresponding projections, $P^1$ and $P^2$.
Then it is well-known (see for instance \cite{KLMV}) that since
\[
P^1 A_k P^2 = P^2 A_k^{[*]} P^1 = 0,
\]
the operators $A_1$ and $A_2$ have the following decompositions
\[
A_k =
\begin{pmatrix}
  A_k^1 & 0
  \\
  i \Phi^{2[*]} \sigma_k \Phi^1 &  A_{k}^2
   \end{pmatrix},
   \qquad
   k=1,2.
\]
Here
$A^1_k = P^1 A_k \restrict{\mathcal P^1}$, $A^2_k = A_k \restrict{\mathcal P^2}$ and
$\Phi^m = \Phi \restrict{\mathcal P^m}$ for $m,k=1,2$.
The first question, is whether the new two collections
\[
\mathcal V^m =
\left(
 \, A^m_1 \, , \,  A^m_2  \, ;  \, \mathcal P^m \, , \,  \Phi^m \, ,  \, E  \, ;  \, \sigma_1 \, ,  \, \sigma_2 \,  , \,  \gamma^m  \, , \, \widetilde{\gamma}^m  \,
\right)
\qquad
m=1,2,
\]
remain vessels.
We note that applying $P_1$ on the right, the input vessel condition becomes
\[
\sigma_1 \Phi^1 A_2^{1[*]}-\sigma_2 \Phi^1 A_1^{1[*]} = \gamma^1 \Phi^1.
\]
This is just the input vessel condition of $\mathcal V_1$ and leads to $\gamma^1 = \gamma$. We are forced to set
\[
\widetilde{\gamma^1} = \gamma + i ( \sigma_1 \Phi^1 \Phi^{1[*]}\sigma_2 - \sigma_2 \Phi^1 \Phi^{1[*]}\sigma_1).
\]
Similarly, applying $P_2$ on the right of the output vessel condition leads to
\[
\sigma_1 \Phi^2 A_2^{2}-\sigma_2 \Phi^2 A_1^{2} = \widetilde{\gamma}^2 \Phi^2
\]
which is the output vessel condition of $\mathcal V_2$.
This implies that $\widetilde{\gamma}^2 = \widetilde{\gamma}$ and furthermore we set
\[
\gamma^2 = \widetilde{\gamma} -  i ( \sigma_1 \Phi^2 \Phi^{2[*]}\sigma_2 -\sigma_2 \Phi^2 \Phi^{2[*]}\sigma_1).
\]
Thus, we have shown the following.
\begin{theorem}[Decomposition theorem]
\label{vesselDecom}
Let $\mathcal V$ be a vessel over a Pontryagin space and let $\mathcal P_2 \subseteq \mathcal P$ be a closed non--degenerate subspace.
Then the two collections
\begin{align*}
\mathcal V_1 & =
\left(\,
P^1 A_1 \restrict{\mathcal P^1}\, , \, P^1 A_2\restrict{\mathcal P^1} \, ;  \,
\mathcal P \ominus \mathcal P_2 \, , \,  \Phi \restrict{\mathcal P \ominus \mathcal P_2}  \, ,  \, E  \, ;  \,
\sigma_1 \, ,  \, \sigma_2 \,  , \,  \gamma^1  \, , \, \widetilde{\gamma}^1  \,
\right)
\\
\mathcal V_2 & =
\left(\,
A_1 \restrict{\mathcal P^2} \, , \,  A_2 \restrict{\mathcal P^2} \, ;  \,
\mathcal P_2 \, , \,  \Phi \restrict{\mathcal P_2}  \, ,  \, E  \, ;  \,
\sigma_1 \, ,  \, \sigma_2 \,  , \,  \gamma^2  \, , \, \widetilde{\gamma}^2  \,
\right)
\end{align*}
are vessels over Pontryagin spaces if
$\gamma^1 = \gamma$, $\widetilde{\gamma}^2 = \widetilde{\gamma}$ and $\widetilde{\gamma}^1 = \widetilde{\gamma}^2$.
\end{theorem}
\begin{pf}
It remain to show that $\widetilde{\gamma}^1 = \widetilde{\gamma}^2$.
The linkage condition of $\mathcal V$ is
\begin{align*}
\widetilde{\gamma} & = \gamma + i( \sigma_1 \Phi \Phi^{[*]} \sigma_2 - \sigma_2 \Phi \Phi^{[*]} \sigma_1)
\\
&
=
\gamma_1 + i( \sigma_1 \Phi_1 \Phi^{[*]}_1 \sigma_2 - \sigma_2 \Phi_1 \Phi^{[*]}_1 \sigma_1)
+ i( \sigma_1 \Phi_2 \Phi^{[*]}_2 \sigma_2 - \sigma_2 \Phi_2 \Phi^{[*]}_2 \sigma_1)
\\
&
=
\gamma_1
+
(
\widetilde{\gamma}_1
-
\gamma_1
)
+
(
\widetilde{\gamma}_2
-
\gamma_2
)
\\
&
=
\widetilde{\gamma}_1
+
\widetilde{\gamma}
-
\gamma_2
\end{align*}
and hence  $\gamma_2 = \widetilde{\gamma}_1$.
\end{pf}

The previous observation leads us to the definition of pairing of two vessels (or equivalently, $2D$ systems)
and determine whether and under which conditions this pairing remains a vessel.
This result is characterized in the following theorem and
the proof in the positive case can be found in \cite{kravitsky1983regular} and in \cite{MR595748}.
\begin{theorem}[Pairing theorem]
Let $\mathcal V_1$ and $\mathcal V_2$ be two commutative vessels with a common external part $(E, \sigma_1, \sigma_2)$, that is,
\[
\mathcal V_k =
\left(\, A^k_1 \, , \,  A^k_2 \, ;  \, \mathcal P^k \, , \,  \Phi^k  \, ,  \, E  \, ;  \,
\sigma_1 \, ,  \, \sigma_2 \,  , \,  \gamma^k  \, , \, \widetilde{\gamma}^k  \,
\right)
\qquad
k=1,2,
\]
where $\mathcal P^1$ is of negative index $\kappa_1$ and $\mathcal P^2$ is of negative index $\kappa_2$.
Then
\[
\mathcal V =
\left(\,
A_1 \, , \,
 A_2 \, ;  \,
\mathcal P^1 \oplus \mathcal P^2 \, , \,  (\Phi^1 \, \, \Phi^2)^T  \, ,  \, E  \, ;  \,
\sigma_1 \, ,  \, \sigma_2 \,  , \,  \gamma  \, , \, \widetilde{\gamma}  \,
\right)
\]
where
\[
A_k  = \begin{pmatrix}
  A^k_1 & 0
  \\
  i \Phi^{2[*]} \sigma_k \Phi^1 &  A^{k}_2
   \end{pmatrix}
   , \qquad
   k=1,2,
\]
is a commutative two-operator vessel over the Pontryagin space $\mathcal P = \mathcal P^1 \oplus \mathcal P^2$ with negative index
$\kappa = \kappa^1 + \kappa^2$
if and only if
$\gamma = \gamma^1$, $\widetilde{\gamma} = \widetilde{\gamma}^2$ and $\gamma^2 = \widetilde{\gamma}^1$.
\end{theorem}
\begin{pf}
Commutativity and the colligation conditions follow from commutativity and the colligation conditions of $\mathcal V_1$ and $\mathcal V_2$.
Then using the following equalities of a vessel
\[
\frac{1}{i}(A_1 A_2^{[*]} - A_2 A_1^{[*]}) = \Phi^{[*]} \gamma \Phi
\qquad {\rm and} \qquad
\frac{1}{i}(A_2^{[*]} A_1 - A_1^{[*]} A_2) = \Phi^{[*]} \widetilde{\gamma} \Phi
,
\]
the sufficient and necessary conditions follow.
\end{pf}

The next result is due to the existence of invariant subspace with maximal negative index of commuting unitary operators \cite{bognar}.
\begin{Cy}
Let $\mathcal V$ be a two-operator commutative vessel over Pontryagin space (with index $\kappa$) such that $A_1$ and $A_2$ are unitary operators.
Then there exists a two-operator commutative vessel over Hilbert space and two-operator commutative vessel over anti--Hilbert space such that $\mathcal V$ is their coupling.
\end{Cy}
\begin{pf}
Using \cite[Lemma 9.2]{bognar} there exists a subspace
$\mathcal P^2 \subseteq \mathcal P$
with $\kappa$ negative index,
which is invariant under the operators $A_1$ and $A_2$.
It remains to apply the decomposition procedure as introduced in Theorem \ref{vesselDecom}.
\end{pf}
%%%%%%%%%%%%%%%%%%%%%%%%%%%%%%%%%%%%
%									%
%	Colligations					%
%									%
%%%%%%%%%%%%%%%%%%%%%%%%%%%%%%%%%%%%
\section{Realization theory of Pontryagin vessels}
\label{secRealization}
\setcounter{equation}{0}

This section is devoted to present the realizations theorems of the various characteristic functions associated to 
a commutative two-operator vessel over Pontryagin space.
In particular, we are mainly interested in the characterization theorem of the normalized joint characteristic function as 
it sets the foundation toward a non--positive de Branges-Rovnyak theory over compact Riemann surfaces, presented in Section \ref{secFuncModel} below.
%%%%%%%%%%%%%%%%%%%%%%%%%%%%%%%%%%%%
%									%
%	Colligations					%
%									%
%%%%%%%%%%%%%%%%%%%%%%%%%%%%%%%%%%%%
\subsection{The complete characteristic function}
Before heading to the characterization theorem of the complete characteristic function, Theorem \ref{JCF_classification},
we must recall the following definitions (we follow here the conventions of \cite{KLMV}).

\begin{definition}
An $n \times n$ matrix function $W(\xi_1,\xi_2,z)$ holomorphic in an open set $K \subseteq \mathbb C ^3$,
is called {\rm consonant} with the determinantal representation $y_1 \sigma_2 - y_2 \sigma_1 + \gamma$
if for any point $(y_1,y_2)$ on the algebraic curve $C$ such that $(\xi_1,\xi_2,\xi_1 y_1 + \xi_2 y_2) \in K$,
the restriction $W(\xi_1,\xi_2,\xi_1 y_1 + \xi_2 y_2)\restrict{E(y)}$ is independent of $\xi_1$ and $\xi_2$.

In such a case, the function $\widehat{W}(y) = W(\xi_1,\xi_2,\xi_1 y_1 + \xi_2 y_2)\restrict{E(y)}$
is called the {\rm trunk} of $W$ with respect to $y_1 \sigma_2 - y_2 \sigma_1 + \gamma$.
\end{definition}

The description of the class of the complete characteristic functions
of vessels over Pontryagin spaces is given in the following theorem.

\begin{theorem}
\label{JCF_classification}
Let $W(\xi_1,\xi_2,z)$ be an $n \times n$ matrix function holomorphic on $\Omega \subset \mathbb C^3$.
Then $W(\xi_1,\xi_2,z)$ is the complete characteristic function of a commutative two-operator vessel 
with a discriminant curve $C$
and input and output determinantal representations $z_1 \sigma_2 - z_2 \sigma_1 + \gamma$ and $z_1 \sigma_2 - z_2 \sigma_1 + \widetilde{\gamma}$,
if and only if:
\begin{enumerate}
\item
\label{JCF_classificationItem1}
$\Omega$ contains 
$\Omega_a = \{(\xi_1,\xi_2,z) \, | \,  |z|^2 > a (|\xi_1|^2 + |\xi_2|^2) \}$  for some $a>0$.
Furthermore, $W(\xi_1,\xi_2,z)$ can be written as
        \begin{equation}
            \label{eqR1}
    W(\xi_1,\xi_2,z) = I -iR(\xi_1,\xi_2,z)(\xi_1 \sigma_1 + \xi_2 \sigma_2),
        \end{equation}
for some $R(\xi_1,\xi_2,z)$ holomorphic on $\Omega$.
\item
\label{JCF_classificationItem2}
For each $\xi_1,\xi_2 \in \mathbb R$ and $z,w \in \mathbb C_+$ such that $(\xi_1,\xi_2,z)\in \Omega$, the kernel
\begin{equation}
\label{kerKappaNeg}
    \frac{W(\xi_1,\xi_2,w)^{[*]}(\xi_1 \sigma_1 + \xi_2 \sigma_2 )W(\xi_1,\xi_2,z)  - (\xi_1 \sigma_1 + \xi_2 \sigma_2)}{-i ( z - \overline{w}) }
\end{equation}
has $\kappa$ negative squares and for $x \in \mathbb C _ {\mathbb R}$ such that $(\xi_1,\xi_2,x)\in \Omega$ 
\begin{equation}
\label{kerKappaNeg3}
    W(\xi_1,\xi_2,x)^{[*]}(\xi_1 \sigma_1 + \xi_2 \sigma_2 ) W(\xi_1,\xi_2,x) = \xi_1 \sigma_1 + \xi_2 \sigma_2.
\end{equation}
\item
\label{JCF_classificationItem3}
$W(\xi_1,\xi_2,z) $ is consonant with $z_1 \sigma_2 - z_2 \sigma_1 + \gamma$ and the trunk $\widehat{W}(y)$
maps $E(y)$ onto $\widetilde{E}(y)$, for any affine $y$ in $C$.
\end{enumerate}
\end{theorem}

    \begin{Rk}
    We may assume, without loss of generality, that $\det \sigma_1 \neq 0$ and $\det \sigma_2 \neq 0$.
    Indeed, since by assumption $C$ does not contain the line in infinity,
    there exist $\xi^{(1)},\xi^{(2)} \in \mathbb R^2$ such that
$\det (\xi^{(1)}_1 \sigma_1 +\xi^{(1)}_2 \sigma_2) \neq 0$,
$\det (\xi^{(2)}_1 \sigma_1 +\xi^{(2)}_2 \sigma_2) \neq 0$
and where
\[
\det
\begin{pmatrix}
\xi^{(1)}_1 & \xi^{(1)}_2\\
\xi^{(2)}_1 & \xi^{(2)}_2
\end{pmatrix}
= 1.
\]
Then it remains to apply the $\alpha$-transformation.
    \end{Rk}
Before presenting the proof, we recall two results which are required in the arguments.
Both hold also in the indefinite inner product case and we refer to \cite[Section 10]{KLMV} for the proofs.

\begin{lem}
Let $R(\xi_1,\xi_2,z)$ be an $n \times n$ matrix function satisfying conditions
\eqref{JCF_classificationItem1}-\eqref{JCF_classificationItem3} in Theorem \ref{JCF_classification}
and let us define
\begin{equation}
\label{eqV}
V(\xi_1,\xi_2) \defEq
- \frac{1}{2 \pi i}
\int_{|z| = r}{ e^{i z} R(\xi_1,\xi_2,z) dz}
,
\end{equation}
where $r >\sup_{\abs{\xi_1}^2+\abs{\xi_2}^2 = 1} \norm{ \xi_1 A_1 + \xi_2 A_2 }  \sqrt{(|\xi_1|^2 + |\xi_2|^2}$.
Then $V(\xi_1,\xi_2)$ satisfies the PDEs
\begin{align}
\label{vFuncPde1}
\frac{\partial V}{\partial \xi_1}\sigma_2 -
\frac{\partial V}{\partial \xi_2}\sigma_1 -
i V \gamma & = 0
\\
\label{vFuncPde2}
\sigma_2 \frac{\partial V}{\partial \xi_1} -
\sigma_1 \frac{\partial V}{\partial \xi_2} -
i \widetilde{\gamma}V & = 0,
\end{align}
and the initial condition
\begin{equation}
\label{vFuncPde3}
i \left(
\sigma_2 V(0,0)\sigma_1
-
\sigma_1 V(0,0)\sigma_2
\right)
=
\widetilde{\gamma} -\gamma
.
\end{equation}
\end{lem}

The next result, needed later in this section, stating that the colligation condition associated to a linear operator $A_2$
is obtained automatically from the colligation condition (of $A_1$) and the input, output and linkage conditions.
\begin{lem}
\label{collCond2ndOper}
Let
$\left( A_2 \, ; \, \mathcal P_2 , \, \Phi_2 \, , \mathbb C ^n \, ; \, \sigma_2 \, \right)$ 
be an irreducible Pontryagin colligation where $\det \sigma_2 \neq 0$.
Let $A_1$ be a linear operator defined on $\mathcal P_2$.
Furthermore, we assume that $A_1$ and $A_2$ commute and that there exist 
selfadjoint operators $\sigma_1, \gamma$ and $\widetilde {\gamma}$ in $E$ such that
\begin{align*}
\Phi ^{[*]} \gamma &= A_2 \Phi ^{[*]} \sigma_1 - A_1 \Phi ^{[*]} \sigma_2  \\
\widetilde{\gamma} \Phi &=\sigma_1 \Phi A_2 - \sigma_2 \Phi A_1  \\
\widetilde{\gamma} &= \gamma + i
\left(
\sigma_1 \Phi \Phi ^{[*]} \sigma_2 - \sigma_2 \Phi \Phi ^{[*]} \sigma_1
\right)
.
\end{align*}
Then
\begin{equation*}
\innerProductIndefiniteReg{A_1 g}{h}
-
\innerProductIndefiniteReg{ g}{A_1 h}
=
i \innerProductIndefiniteReg{\sigma_1 \Phi g}{\Phi h }
\end{equation*}
holds for every $g,h \in \mathcal P_2$.
\end{lem}
\begin{pf}
Noting that the colligation is, by assumption, irreducible,
the proof is similar to the proof in the Hilbert space setting, see \cite[Proposition 10.4.8]{KLMV}.
\end{pf}

\begin{pf}[of Theorem \ref{JCF_classification}]
We start with the sufficient part.
Let $\mathcal V$ be a vessel over a Pontryagin space $\mathcal P$.
Then the complete characteristic function of $\mathcal V$ is analytic in
\[
\Omega_d  = \left \{(\xi_1,\xi_2,z) \, \big| \, |z| > d \sqrt{\abs{\xi_1}^2+\abs{\xi_2}^2}\right \}
\]
where $d=\sup_{\abs{\xi_1}^2+\abs{\xi_2}^2 = 1} \norm{ \xi_1 A_1 + \xi_2 A_2 }$.
Furthermore, we note that the collection
\[
\mathcal C =
\left( \xi_1 A_1 +\xi_2 A_2 \, ; \, \mathcal P , \, \Phi \, , E \, ; \, \xi_1 \sigma_1 + \xi_2 \sigma_2 \, \right),
\]
is a single-operator colligation over the Pontryagin space $\mathcal P$, with characteristic function 
    $W(\xi_1 \sigma_1 + \xi_2 \sigma_2,z)$ where $\det {(\xi_1 \sigma_1 + \xi_2 \sigma_2) } \neq 0$.
Thus, applying Theorem \ref{singleOpColl}, the kernel \eqref{kerKappaNeg}
has at most $\kappa$ negative squares and for some $\xi_1$ and $\xi_2$ has exactly $\kappa$ negative squares.
Lastly, part \eqref{JCF_classificationItem3} is followed by Lemma \ref{jcfLem}.
This completes the proof of the "only if" part.
\smallskip

We continue with the necessary part. 
For the sake of clarity, we divide this part of the proof into a number of steps.

\setcounter{Step}{0}
\begin{Step}
Construction of a single-operator colligation over Pontryagin space.
\end{Step}
By assumption, in the direction $\xi = (0,1)$, the kernel
\[
    \frac{W(0,1,w)^{[*]}  \sigma_2 W(0,1,z)  -  \sigma_2}{-i ( z - \overline{w}) }
\]
has $\kappa$ negative squares.
Hence, using Theorem \ref{singleOpColl}, there exists an irreducible single-operator colligation over Pontryagin space
\[
\mathcal C _2 = \left( A_2 \, ; \, \mathcal P_2 , \, \Phi_2 \, , \mathbb C ^n \, ; \, \sigma_2 \, \right),
\]
where $\mathcal P_2$ is a Pontryagin space of negative index $\kappa$,
$A_2$ is bounded in $\mathcal P_2$ and
$\Phi_2$ is a bounded operator from $\mathcal P_2$ to $\mathbb C ^n$.
The characteristic function of the colligation $\mathcal C _2$ is $W(0,1,z)$
and thus we have
\[
W(0, \xi , z ) = I - i \Phi_2 (\xi A_2 - zI)^{-1} \Phi_2^{[*]} \xi \sigma_2
.
\]
Furthermore, we note that \eqref{eqV} satisfies
\begin{align}
\label{coll2V}
V(0, \xi )
= &
- \frac{1}{2 \pi i} \int_{|z| = r}{ e^{i z} \Phi_2 (\xi A_2 - zI)^{-1} \Phi_2^{[*]} dz}
\\ = &
\Phi_2 \exp(i \xi A_2) \Phi_2^{[*]}
.
\nonumber
\end{align}

\begin{Step}
Study the space of solutions of the PDE
\begin{equation}
\label{vPde1}
\sigma_2 \frac{\partial v}{\partial t_1}
- \sigma_1 \frac{\partial v}{\partial t_2}
+ i \widetilde{\gamma} v = 0
,
\end{equation}
with the initial condition
    \begin{equation}
        \label{eqInitCond}
v_p(0,t_2) = -i \Phi_2 \exp {(it_2 A_2)}p.
    \end{equation}
\end{Step}
For each $p \in \mathcal P_2$, we associate a unique solution which is the restriction to $\mathbb R^2$ of an entire funtion, 
denoted by $v_p(t_1,t_2)$, to the PDE \eqref{vPde1} with the initial condition \eqref{eqInitCond}.
First, we consider $p = A_2^k \Phi_2^{[*]}e$ for some $e \in \mathbb C ^n$ and $k \in \mathbb N$.
Then, using \eqref{vFuncPde2}, the (unique) solution is given by
\begin{equation}
\label{vSolPde}
v_p(t_1,t_2) = (-i
)^{k+1} \frac{\partial ^k}{\partial t_2^k} V(t_1,t_2)e.
\end{equation}
The initial condition \eqref{eqInitCond} follows from \eqref{coll2V}.
For uniqueness, let assume that $v_p(0,t_2)=0$ for all $t_2 \in \mathbb R$.
Then by \eqref{vPde1} and since $\det{\sigma_2} \neq 0$, we can conclude that $\frac{\partial}{\partial t_1}V(0,t_2) =0$.
By induction on $k_1$ and $k_2$ we have $\frac{\partial^{k_1+k_2}}{\partial  t_1 ^{k_1} \partial t_1^{k_2}}V(0,t) = 0$.
By assumption, $v(z_1,z_2)$ is entire and hence must be zero.
\smallskip

By linearity, one can extend the statement and associate a solution also for every $p$ in the linear envelope $A_2^k \Phi_2^{[*]} e$.
It remains to take a sequence $(p_n)_{n=1}^\infty$ in the linear envelope $A_2^k \Phi_2^{[*]} e$ converging to some $p$ in $\mathcal P_2$, 
and show that the limit of the sequence $(v_{p_n})_{n=1}^\infty$, a sequence of solutions for \eqref{vPde1}, exists.
Noting that for any $t_2 \in \mathbb R$
\[
    v_{p_n}(0,t_2) = 
    -i \Phi_2 \exp {(it_2 A_2)}p_n
    \xrightarrow{n \rightarrow \infty}
    -i \Phi_2 \exp {(it_2 A_2)}p
.
\]
Then the limit indeed exists by \cite[Proposition 10.4.7]{KLMV}
and noting that the proposition is given in terms of the norms of elements of the sequence and not the indefinite inner product.

\begin{Step}
Construction of an isometry between $\mathcal P_2$ and the space of solutions of \eqref{vPde1}.
\end{Step}
Let us define the space of solutions of the PDE in \eqref{vPde1} by
\[
\mathcal P = \left \{ \nu_p(t_1,t_2) \, | \, p \in \mathcal P_2 \right \}
.
\]
By Lemma \ref{lemmaPrincSubSpace} \eqref{princSubSpace2}, the mapping
$U_2 : \mathcal P_2 \rightarrow \mathcal P$ given by $U_2 (h) = \nu_h(t_1,t_2)$, is injective.
Hence $\mathcal P$ becomes a Pontryagin space (with $\kappa$ negative squares) with the indefinite inner product inherited from $\mathcal P_2$ 
and furthermore, this identification makes $U_2$ to an isometry.
\smallskip

We note that under the isometry $U_2$, the operators $\mathcal A_2, \Psi$ and $\Psi ^{[*]}$ can be written as
\begin{align}
\label{eqCol1}
\Psi \nu_p & = \Phi_2 p  = i \nu _p (0,0) \\
\Psi^{[*]} e & = \nu_{\Phi^{[*]}_2e}  = - i V(t_1,t_2)e
\label{eqCol2}
\end{align}
and
\[
\mathcal A_2 \nu_p  = \nu_{A_2 p}  =  -i\frac{\partial}{\partial t_2} \nu_p,
    \]
respectively.
The collection
\[
\mathcal C _2
=
\left( \,
\mathcal A_2 = -i\frac{\partial}{\partial t_2} \, ;
\, \mathcal P 	\, ,
\, \Psi 	\, ,
\, \mathbb C ^n \, ;
\, \sigma_2  \,
\right)
\]
forms an irreducible colligation over Pontryagin space $\mathcal P$ with $\kappa$ negative index and with the characteristic function $W(0,1,z)$.

\begin{Step}
The operator $\mathcal A_1 \defEq -i\frac{\partial}{\partial t_1} $ is a linear operator in $\mathcal P$ and defined on the linear span of functions in
$ \bigvee_{n=0}^{\infty}{\mathcal A_2^k \Psi^{[*]} e}$.
\end{Step}
Using \eqref{vFuncPde1}
\[
\frac{\partial}{\partial t_1} V(t_1,t_2) \sigma_2 -
\frac{\partial}{\partial t_2} V(t_1,t_2) \sigma_1 +
i V(t_1,t_2) \gamma = 0,
\]
we have, since $\sigma_2$ is invertible, the following
\[
\frac{\partial^{k+1}}{\partial t_1 \partial t_2^{k}} V(t_1,t_2)  =
\frac{\partial^{k+1}}{\partial t_2^{k+1}} V(t_1,t_2) \sigma_1 \sigma_2^{-1}
- i \frac{\partial^{k}}{\partial t_2 ^{k}} V(t_1,t_2) \gamma \sigma_2^{-1}.
\]
Hence $\mathcal A_1$ becomes a linear operator on $\mathcal P$ defined on the linear span of $\mathcal A_2 ^k \Psi^{[*]}e$.

\begin{Step}
The collection
\begin{equation}
\label{shiftVessel}
\mathcal V
=
\left(
-i\frac{\partial}{\partial t_1} \, , \, -i\frac{\partial}{\partial t_2} \, ;
\, \mathcal P 	\, , \, \Psi 	\, , \, \mathbb C ^n \, ;
\, \sigma_1 	\, , \, \sigma_2 	\, , \,	\gamma \, , \, \widetilde{\gamma} \,
\right)
\end{equation}
is a commutative two--operator Pontryagin vessel with the complete characteristic function $W(\xi_1,\xi_2,z)$.
\end{Step}
\smallskip

First, we note that $-i\frac{\partial}{\partial t_1}$ and $-i\frac{\partial}{\partial t_2}$ commute.
Furthermore, the input vessel condition
\begin{equation}
\label{shiftVesselInput}
\sigma_1 \Psi \left(-i\frac{\partial}{\partial t_2} \right) -
\sigma_2 \Psi \left(-i\frac{\partial}{\partial t_1} \right)
=
\widetilde{\gamma} \Psi
\end{equation}
is exactly \eqref{vPde1} evaluated at $t = (0,0)$.
The output condition 
\begin{equation}
\label{shiftVesselOutput}
\left(-i\frac{\partial}{\partial t_2} \right) \Psi^{[*]} \sigma_1   -
\left(-i\frac{\partial}{\partial t_1} \right) \Psi^{[*]} \sigma_2
=
\Psi ^{[*]} \gamma
\end{equation}
is follows by \eqref{vFuncPde2} and 
the linkage vessel condition
\begin{equation}
\label{shiftVesselLinkage}
\widetilde{\gamma}  - \gamma
=
i \left( \sigma_1 \Psi\Psi^{[*]} \sigma_2   - \sigma_2 \Psi\Psi^{[*]} \sigma_1 \right),
\end{equation}
is exactly \eqref{vFuncPde3}.
Combining the colligation condition corresponding to $\mathcal A_2$ with \eqref{shiftVesselInput}, 
\eqref{shiftVesselOutput} and \eqref{shiftVesselLinkage},
we apply Lemma \ref{collCond2ndOper} and conclude the colligation condition corresponding to $\mathcal A_1$, i.e.
\[
i
\innerProductIndefiniteTri{\sigma_1 \Psi \nu_1}{\Psi \nu_2}{\mathbb C ^n}
=
\innerProductIndefiniteTri{\mathcal A_1 \nu_1}{\nu_2}{\mathcal P}-
\innerProductIndefiniteTri{\nu_1}{\mathcal A_1 \nu_2}{\mathcal P}
\]
holds for any $\nu_1,\nu_2 \in \mathcal P$.
\begin{Step}
The operator $\mathcal A_1$ is bounded in $\mathcal P$.
\end{Step}
To show that the shift-vessel $\mathcal V$ introduced in \eqref{shiftVessel} is indeed a vessel,
it remains to show that the linear operator $\mathcal A_1$ is bounded.
Applying Theorem \ref{singleOpColl} once again, but now with respect to the function $W(1,0,z)$,
we conclude that there exists an irreducible single-operator colligation over Pontryagin space
\[
\mathcal C _1 = \left( A_1 \, ; \, \mathcal P_1 , \, \Phi_1 \, , \mathbb C ^n \, ; \, \sigma_1 \, \right),
\]
with the characteristic function $W(1,0,z)$.
Hence
\begin{align*}
W(\xi,0,z) & = I - i \Phi_1 (\xi A_1 - zI)^{-1} \Phi_1^{[*]} \xi \sigma_1,
\end{align*}
and we also have
\begin{align}
\label{eqVc1}
V(\xi,0)
= & 
- \frac{1}{2 \pi i} \int_{|z| = r}{ e^{i z} \Phi_1 (\xi A_1 - zI)^{-1} \Phi_1^{[*]} dz}
\\  = &
\Phi_1 \exp{i \xi A_1} \Phi_1^{[*]}
.
\nonumber
\end{align}
Thus, using \eqref{eqVc1}, \eqref{eqCol1} and \eqref{eqCol2}, we have the following identity
\[
\Psi \mathcal A_1 ^k \Psi^{[*]} e = i^k \frac{\partial^k V}{\partial t_1^k}(0,0)e = \Phi_1 A_1 ^k \Phi_1 ^{[*]} e.
\]
By induction on $l$ and using the colligation conditions associated to $\mathcal C _1$, we have
\[
\innerProductIndefiniteReg{ \mathcal A_1 ^k \Psi^{[*]} e_1 }{\mathcal A_1 ^l \Psi^{[*]} e_2}
=
\innerProductIndefiniteReg{ A_1 ^k \Phi^{[*]} e_1 }{A_1 ^l \Phi^{[*]} e_2} .
\]
Thus, the mapping
\[
U_1 \, : \, A_1^k \Phi^{[*]} e \rightarrow \mathcal A_1^k \Psi^{[*]} e
\]
defines an isometry from the linear envelop of $A_1^{k} \Phi_1^{[*]} e$ into $\mathcal P$.
Since $\mathcal C_1$ is a irreducible colligation, the isometry is defined on a dense subset of $\mathcal P_1$.
Furthermore, since the solutions of \eqref{vPde1} are restrictions to $\mathbb R^2$ of entire functions, 
the linear envelope $\mathcal A_1 ^l \Psi^{[*]} e$ is dense in $\mathcal P$, which also has $\kappa$ negative squares.
Then, and this is a crucial point, based on the isometry extension theorem in the Pontryagin space setting 
(see \cite[Theorem 4.1.2]{adrs}, in a more general setting)
there exists a continuous extension of $U_1$ to $\mathcal P$.
Hence $\mathcal A_1$ can be extended to a bounded linear operator on $\mathcal P$.
\begin{Step}
The complete characteristic function of the vessel $\mathcal V$ is equal to $W(\xi_1,\xi_2,z)$.
\end{Step}
Let $W_{\mathcal V}(\xi_1,\xi_2,z)$ be the complete characteristic function of the vessel \eqref{shiftVessel}.
Then $W_{\mathcal V}(\xi_1,\xi_2,z)$ and $W(\xi_1,\xi_2,z)$ coincide along $\xi = (0,1)$.
Using that the complete characteristic function 
is fully determined by its value for any fixed values of $\xi_1,\xi_2$, see Corollary \ref{sUniqTrunk} below,
the proof is completed.
\end{pf}

%%%%%%%%%%%%%%%%%%%%%%%%%%%%%%%%%%%%%
%									%
%	JCF								%
%									%
%%%%%%%%%%%%%%%%%%%%%%%%%%%%%%%%%%%%%
\subsection{Realization of the joint characteristic function}
In this section we answer the question whether and under which assumptions a mapping between two vector bundles 
associated to two determinantal representations defined on the same curve, is the joint characteristic function of a Pontryagin space vessel.
The main tool that we use throughout the section is referred as the restoration formula.
The restoration formula allows us to recover entirely the complete characteristic function at a point $z$ by
the evaluations of the joint characteristic function in $n$ distinct points intersecting the line $z = \xi_1 \lambda_1 + \xi_2 \lambda_1$ and the curve $C$.
Let $(\xi_1,\xi_2,z)$ be an arbitrary point in $\mathbb C ^3$
and we consider the $n$ distinct points $(\lambda_1^{(j)},\lambda_2^{(j)})\in \mathbb R ^2$ such that
$\xi_1 \lambda_1^{(j)} + \xi_2 \lambda_1^{(j)} = z$.
The eigenspaces $E (\lambda^{(j)})$, assuming maximality, satisfy
\[
E (\lambda^{(1)}) + \cdots + E (\lambda^{(n)}) = E,
\]
and together the corresponding projections, denoted by $P(\lambda^{(j)};\xi_1,\xi_2,z)$,
we have the following ({\it restoration formula}) decomposition 
\[
W(\xi_1,\xi_2,z) =
\sum_{\xi_1 \lambda^{(j)}_1 + \xi_2 \lambda^{(j)}_2 = z}^{}{S(\lambda^{(j)})P(\lambda^{(j)};\xi_1,\xi_2,z)}.
\]
The upcoming lemma illustrates an important consequence of the restoration formula.
\begin{lem}
\label{sUniqTrunk}
Let $W(\xi_1,\xi_2,z)$ be an $n \times n$ matrix function holomorphic on an open set $K \subset \mathbb C ^3$ 
and consonant with a selfadjoint determinantal representation $\lambda_1 \sigma_2 - \lambda_2 \sigma_1 +\gamma$ 
of a real projectiveu reduced plane curve $C$ of degree $n$ which does not contain the line at infinity.
Then, for a fixed $\xi_1^0,\xi_2^0 \in \mathbb R$ such that $\det(\xi_1^0 \sigma_1 + \xi_2^0 \sigma_2) \neq 0$,
the function $W(\xi_1,\xi_2,z)$ is uniquely determined by the values $W(\xi^0_1,\xi^0_2,z)$. 
% Hence $W(\xi_1,\xi_2,z)$ is uniquely determined by its trunk corresponding to
% $\lambda_1 \sigma_2 - \lambda_2 \sigma_1 +\gamma$.
\end{lem}

The proofs of the restoration formula and Lemma \ref{sUniqTrunk} are similar to the classical case.
\smallskip

The metric properties \eqref{kerKappaNeg} and \eqref{kerKappaNeg3} of the complete characteristic function
should be translated to corresponding metric properties of the joint characteristic function.
The metric is given in term of a Hermitian pairing between fibers of kernel bundles corresponding to non--conjugate points
\[
\innerProductIndefiniteQuad{u}{v}{\lambda_1,\lambda_2}{E} =
i \frac{v^*(\xi_1 \sigma_1 + \xi_2 \sigma_2 )u}{\xi_1(\lambda_1^1 - \overline{\lambda}_1^2)  + \xi_2(\lambda_2^1 - \overline{\lambda}_2^2) },
\]
where $u \in E (\lambda^1)$ and $v \in E (\lambda^2)$. 
This pairing, in the case of conjugate points, is given by
\[
\innerProductIndefiniteQuad{u}{v}{\lambda_1,\lambda_2}{E} =
i \frac{v^*(\xi_1 \sigma_1 + \xi_2 \sigma_2 )u}{\xi_1 d\lambda_1 (\lambda)  + \xi_2 d\lambda_2(\lambda) }
.
\]
Then, the properties of CCF relative to the indefinite inner product
are equivalent to indefinite scalar properties on $C$.
This point is illustrated in the following lemma.
\begin{theorem}
\label{realizationPreLem}
Let $S(\lambda)$ be a mapping between the vector bundles $E$ and $\widetilde{E}$,
holomorphic on an open set $\Omega \subseteq C$ containing
$\Omega_a = \{ (x_0,x_1,x_2) \in C \, | \, \abs{x_1}^2 +\abs{x_2}^2 > a \abs {x_0}^2\}$ for some $a>0$
and meromorphic on the complement of the real point of $C$.
Let $W(\xi_1,\xi_2,z)$ be an $n \times n$ matrix function holomorphic on $K \subseteq \mathbb C^3$
containing $K_a = \{ (\xi_1,\xi_2,z) \in \mathbb C^3 \, | \, a \left( \abs{\xi_1}^2 +\abs{\xi_2}^2 \right) < \abs {z}^2\}$,
which is consonant with $\lambda_1 \sigma_2 - \lambda_2 \sigma_1 + \gamma$ and $W(\xi_1,\xi_2,\xi_1 \lambda_1 + \xi_2 \lambda_2)$ coincides with $S$
and meromorphic for $z$ on the complement of the real axis.
Then for $N>0$ and any choice of $\lambda^{(1)},...,\lambda^{(N)} \in \Omega$
\begin{equation}
\label{jcfCond1}
\left(
\innerProductIndefiniteQuad{S(\lambda^{(h)})u^{(h)}}{S(\lambda^{(j)})u^{(j)}}{\lambda^{(h)},\lambda^{(j)}}{\widetilde{E}}
-
\innerProductIndefiniteQuad{u^{(h)}}{u^{(j)}}{\lambda^{(h)},\lambda^{(j)}}{E}
\right)_{h,j=1}^{N}
\end{equation}
has at most $\kappa$ negative eigenvalues where $(u^{(h)}\in E(\lambda^{(h)}); h=1,...,N)$, and
\begin{equation}
\label{jcfCond2}
\innerProductIndefiniteReg{S(\lambda)u}{S(\overline{\lambda})v}^{\widetilde{E}}_{\lambda,\overline{\lambda}}
=
\innerProductIndefiniteReg{ u }{ v }^{E}_{\lambda,\overline{\lambda}}
\end{equation}
for $(u \in E(\lambda), v \in E(\overline{\lambda}))$,
if and only if
for each $\xi_1,\xi_2 \in \mathbb R$ and for all $w,z \in \mathbb C_+$ such that $(\xi_1,\xi_2,w),(\xi_1,\xi_2,z) \in K$
the kernel \eqref{kerKappaNeg}
has $\kappa$ negative squares and \eqref{kerKappaNeg3} holds.
\end{theorem}
\begin{pf}
We fix $\xi_1, \xi_2 \in \mathbb R$ and consider $z \in \mathbb C \setminus \mathbb R$ 
such that $(\xi_1, \xi_2 , z) \in K$ and such that the line $\xi_1 \lambda_1 + \xi_2 \lambda_2 = z$ intersects
with the curve $C$ at $n$ distinct points (denoted by $\lambda^{(1)},...,\lambda^{(n)})$.
Then, using the restoration formula, on one hand we have
\begin{align}
    v_k^{[*]}
    W(\xi_1, & \xi_2  , z_k)^{[*]}
    \frac{\xi_1 \sigma_1 +\xi_2 \sigma_2}{-i(z_l - \overline{z_k})}
W(\xi_1, \xi_2 , z_l)
v_l
=
\nonumber
\\
= &
\sum_{h,j=1}^{n}
\left(\widehat{W}(\lambda_k^{(j)})P(\lambda_k^{(j)};\xi,z)v_k \right)^{[*]}
\frac{\xi_1 \sigma_1 + \xi_2 \sigma_2}{-i(z_l - \overline{z_k})}
\times
\nonumber
    \\
&
\times
\left(\widehat{W}(\lambda_l^{(h)})P(\lambda_l^{(h)};\xi,z)v_l \right)
\nonumber
\\
\label{5_8B}
= &
\sum_{h,j=1}^{n}
{
\innerProductIndefiniteQuad
{S(\lambda_l^{(h)})P(\lambda_l^{(h)};\xi,z)v_l}
{S(\lambda_k^{(j)})P(\lambda_k^{(j)};\xi,z)v_k}
{\lambda_l^{(h)},\lambda_k^{(j)}}
{\widetilde{E}}
}.
\end{align}
On the other hand
\begin{align}
v_k^{[*]}
\frac{(\xi_1 \sigma_1 +\xi_2 \sigma_2)}{-i(z_l - \overline{z_k})}
v_l
=
\sum_{h,j=1}^{n}
{
\innerProductIndefiniteQuad
{P(\lambda_l^{(h)};\xi,z)v_l}
{P(\lambda_k^{(j)};\xi,z)v_k}
{\lambda_l^{(h)},\lambda_k^{(j)}}
{E}
}
\label{5_8A}
.
\end{align}
Hence, combining \eqref{5_8B} and \eqref{5_8A} and using that $S$ coincide with the trunk of $W$, we may conclude the following:
\begin{align*}
v_k^{[*]}
&
\frac{W(\xi_1, \xi_2 , z_k)^{[*]}(\xi_1 \sigma_1 +\xi_2 \sigma_2)W(\xi_1, \xi_2 , z_l)- (\xi_1 \sigma_1 +\xi_2 \sigma_2)}{-i(z_l - \overline{z_k})}
v_l =
\\ = &
\sum_{h,j=1}^{n}
{ }
\innerProductIndefiniteQuad
{S(\lambda_l^{(h)})P(\lambda_l^{(h)};\xi,z)v_l}
{S(\lambda_k^{(j)})P(\lambda_k^{(j)};\xi,z)v_k}
{\lambda_l^{(h)},\lambda_k^{(j)}}
{\widetilde{E}}
-
\\
& -
\innerProductIndefiniteQuad
{P(\lambda_l^{(h)};\xi,z)v_l}
{P(\lambda_k^{(j)};\xi,z)v_k}
{\lambda_l^{(h)},\lambda_k^{(j)}}
{E}
.
\end{align*}
Thus, if the matrix \eqref{jcfCond1} has $\kappa$ negative eigenvalues, then the kernel \eqref{kerKappaNeg} also has $\kappa$ negative squares.
\smallskip

We now turn and choose $x \in \mathbb R$
such that there are $n$ distinct affine points, denoted by $(\lambda_1^{(j)},\lambda_2^{(j)})_{j=1}^{n}$, 
satisfying $\xi_1\lambda_1 ^{(j)} + \xi_2 \lambda_2 ^{(j)} = x$.
Then, in particular, for $\lambda^{(j)} \neq \overline{\lambda^{(k)}}$ it follows that
\[
\xi_1 (\lambda^{(j)}_1 - \overline{\lambda^{(k)}_1}) + \xi_2 (\lambda^{(j)}_2 - \overline{\lambda_2^{(k)}}) = 2 \imag x = 0.
\]
As a consequence, 
for any $v \in E(\lambda^{(k)})$ and any $u \in E(\lambda^{(j)})$ 
(and similiarliy for any $v \in \widetilde{E}(\lambda^{(k)})$ and $u \in \widetilde{E}(\lambda^{(j)})$), 
a short computation yields the following
\begin{equation}
\label{eqInnerProZero}
v^{[*]}(\xi_1 \sigma_1 +\xi_2 \sigma_2)u = 0.
\end{equation}
Applying the restoration formula once again, but now on \eqref{jcfCond2}, we have the following:
\begin{align*}
v^{[*]}  W & (\xi_1, \xi_2 , z)^{[*]} (\xi_1 \sigma_1 +\xi_2 \sigma_2)W(\xi_1, \xi_2 , z)v
\\
= &
\sum_{\substack{j,k: \\ \lambda^{(j)} = \overline{\lambda^{(k)}} }}
v^{[*]}
\left[
\widehat{W}(\lambda^{(k)})P(\lambda^{(k)};\xi,z)
\right]^{[*]}
(\xi_1\sigma_1 + \xi_2 \sigma_2)
\widehat{W}(\lambda^{(j)})P(\lambda^{(j)};\xi,z)u.
\end{align*}
Since $S$ coincides with the trunk of $W$ and using \eqref{eqInnerProZero}, we have
\begin{align*}
v^{[*]}  W & (\xi_1, \xi_2 , z)^{[*]} (\xi_1 \sigma_1 +\xi_2 \sigma_2)W(\xi_1, \xi_2 , z)u
\\
= &
\sum_{\substack{j,k: \\ \lambda^{(j)} = \overline{\lambda^{(k)}} }}
\left[ S(\lambda^{(k)})P(\lambda^{(k)};\xi,z) v \right] ^{[*]}
(\xi_1\sigma_1 + \xi_2 \sigma_2)
S(\lambda^{(j)})P(\lambda^{(j)};\xi,z)u
\\
= &
v^{[*]}
(\xi_1\sigma_1 + \xi_2 \sigma_2)
u
,
\end{align*}
and \eqref{kerKappaNeg3} follows.
\smallskip

Conversely, let us assume that \eqref{kerKappaNeg} and \eqref{kerKappaNeg3} hold. 
Then \eqref{kerKappaNeg3} can be extended to every $z \in \mathbb C$ by
\begin{equation}
\label{kerKappaNeg2Z}
W(\xi_1,\xi_2,\overline{z})^{[*]}(\xi_1 \sigma_1 + \xi_2 \sigma_2 ) W(\xi_1,\xi_2,z) = \xi_1 \sigma_1 + \xi_2 \sigma_2.
\end{equation}
Let $\lambda = (\lambda_1,\lambda_2)$ be a non-real affine point on $C$ such that $\lambda, \overline{\lambda} \in \Omega$.
We then choose $\xi_1,\xi_2 \in \mathbb R$ such that
$(\xi_1,\xi_2,\xi_1 \lambda_1 + \xi_2 \lambda_2)$ and $(\xi_1,\xi_2,\xi_1 \overline{\lambda_1} + \xi_2 \overline{\lambda_2})$ are both belong to $K$.
Then for $u \in E(\lambda)$ and $v \in E(\overline{\lambda})$, using \eqref{kerKappaNeg2Z}, we have
\begin{align*}
\innerProductIndefiniteReg{S(\lambda)u}{S(\overline{\lambda})v}^{\widetilde{E}}_{\lambda,\overline{\lambda}}
= &
\innerProductIndefiniteReg{W(\xi_1,\xi_2,\xi_1 \lambda_1 + \xi_2 \lambda_2)u}{W(\xi_1,\xi_2,\xi_1 \overline{\lambda}_1 + \xi_2 \overline{\lambda}_2) v}^{\widetilde{E}}_{\lambda,\overline{\lambda}}
\\
= &
\frac
{v^{[*]} W(\xi,\xi \overline{\lambda})^{[*]} (\xi_1 \sigma_1 +\xi_2 \sigma_2)W(\xi,\xi \lambda) u }
{\xi_1 d \lambda _1 (\lambda) +\xi_2 d\lambda _2 (\lambda)}
\\
= &
\frac
{v^{[*]}  (\xi_1 \sigma_1 +\xi_2 \sigma_2)  u }
{\xi_1 d \lambda_1 (\lambda) +\xi_2 d \lambda _2 (\lambda)}
=
\innerProductIndefiniteReg{u}{v}^{\widetilde{E}}_{\lambda,\overline{\lambda}}
\end{align*}
and \eqref{jcfCond2} follows. \smallskip

We fix $N>0$, then by assumption, for any $\xi_1,\xi_2 \in \mathbb R$, any $z_1,...,z_N \in \mathbb C$ such that
$(\xi_1,\xi_2,z_j) \in K$ and for any $v_1,...,v_N \in \mathbb C^n$, by \eqref{kerKappaNeg} the matrix
\begin{equation}
\label{eqKer12}
\left(
v_j ^{[*]}
\frac
{W(\xi_1,\xi_2,z_j)^{[*]}(\xi_1 \sigma_1 + \xi_2 \sigma_2 ) W(\xi_1,\xi_2,z_k) - (\xi_1 \sigma_1 + \xi_2 \sigma_2)}
{-i ( z - \overline{w}) }
v_k
\right)
_{j,k=1}^{N}
\end{equation}
has at most $\kappa$ negative squares.
\smallskip

Let $\lambda^{(j)}_1,\lambda^{(j)}_2$ be affine representations of points in $\Omega$ 
such that $\lambda^{(j)} \neq \overline{\lambda^{(k)}}$ whenever $k \neq j$.
Then, we can choose $\xi_1,\xi_2 \in  \mathbb R$ such that for all $1 \leq j \neq k \leq N$ we have the following
\[
(\xi_1 , \xi_2 , \xi_1 \lambda^{(j)}_1 + \xi_2 \lambda^{(j)}_2) \in K,
\]
and
\[
\xi_1 (\lambda^{(j)}_1 - \overline{\lambda^{(k)}_1}) + \xi_2 (\lambda^{(j)}_2 - \overline{\lambda^{(k)}_2}) \neq 0.
\]
That is, the selected points belong to the region of analyticity of $W$ and are above distinct points relative to the direction $(\xi_1,\xi_2)$.
Thus we have
\begin{align*}
&
\left(
\innerProductIndefiniteQuad
{S(\lambda^{(j)}) v_j}
{S(\lambda^{(k)}) v_k}
{\lambda^{(k)},\lambda^{(j)}}
{\widetilde{E}}
-
\innerProductIndefiniteQuad
{v_j}
{v_k}
{\lambda^{(k)},\lambda^{(j)}}
{E}
\right)
_{j,k=1}^{N}
\\ = &
\left(
\innerProductIndefiniteQuad
{W(\xi, \xi \lambda^{(j)}) v_j}
{W(\xi, \xi \lambda^{(k)}) v_k}
{\lambda^{(k)},\lambda^{(j)}}
{\widetilde{E}}
-
\innerProductIndefiniteQuad
{v_j}
{v_k}
{\lambda^{(k)},\lambda^{(j)}}
{E}
\right)
_{j,k=1}^{N}
\\ = &
\left(
    v_k ^{[*]}
\frac{
        W(\xi,\xi \lambda^{(k)})^*
        (\xi_1 \sigma_1 + \xi_2 \sigma_2 )
        W(\xi,\xi \lambda^{(j)}) - (\xi_1 \sigma_1 + \xi_2 \sigma_2)
    }
{\xi_1 (\lambda^{(j)}_1 - \overline{\lambda^{(k)}_1}) + \xi_2 (\lambda^{(j)}_2 - \overline{\lambda^{(k)}_2})}
v_j \right)
_{j,k=1}^{N}
.
\end{align*}
and due to \eqref{eqKer12}, the result follows.
\end{pf}

We now can present a theorem which characterizes the family of joint characteristic functions.
\begin{theorem}
\label{realizationTh}
Let $C$ be a plane projective curve of degree $n$ with two selfadjoint representations
$\lambda_1 \sigma_2 - \lambda_2 \sigma_1 + \gamma$
and $\lambda_1 \sigma_2 - \lambda_2 \sigma_1 + \widetilde{\gamma}$,
and two corresponding vector bundles $E(x)$ and $\widetilde{E}(x)$.
A mapping $S$ of the vector bundles $E$ and $\widetilde{E}$ on $C$
is the joint characteristic function of a commutative two-operator Pontryagin space vessel with $\kappa$ negative squares,
with discriminant curve $C$, input determinantal representation $\det(\lambda_1\sigma_2 - \lambda_2 \sigma_1 + \gamma)$
and output determinantal representation $\det(\lambda_1\sigma_2 - \lambda_2 \sigma_1 + \widetilde{\gamma}) $
if and only if:
\begin{enumerate}
\item
% \label{realizationThA}
$S$ equals to identity at infinity and is holomorphic in a neighborhood of the points of $C$ at infinity.
\item
\label{realizationThB}
$S$ is meromorphic on the complement of the set of real points of $C$ and for all affine points 
$\lambda,\lambda^{(1)},...,\lambda^{(N)}$ on $C$ in its region of analyticity 
($\lambda^{(h)} \neq \lambda^{(j)}$, $\overline{\lambda}$ in the region of analyticity of $S$)
\[
\left(
\innerProductIndefiniteQuad{S(\lambda^{(h)})u^{(h)}}{S(\lambda^{(j)})u^{(j)}}{\lambda^{(h)},\lambda^{(j)}}{\widetilde{E}}
-
\innerProductIndefiniteQuad{u^{(h)}}{u^{(j)})}{\lambda^{(h)},\lambda^{(j)}}{E}
\right)_{h,j=1}^{N}
\]
has at most (and for some choice exactly) $\kappa$ negative squares where $(u^{(h)}\in E(\lambda^{(h)}); h=1,...,N)$, and
\[
\innerProductIndefiniteReg{S(\lambda)u}{S(\overline{\lambda})v}^{\widetilde{E}}_{\lambda,\overline{\lambda}}
=
\innerProductIndefiniteReg{ u }{ v }^{E}_{\lambda,\overline{\lambda}}
\]
for $(u \in E(\lambda), v \in E(\overline{\lambda}))$.
\end{enumerate}
\end{theorem}

\begin{pf}
There exists, see \cite[Proposition 10.5.3]{KLMV}, an $n \times n$ matrix function, denoted by $W(\xi_1,\xi_2,z)$,
holomorphic on the open set
\begin{align*}
K=\{ (\xi_1,\xi_2,z)\, | \, & (\xi_1,\xi_2,z) \neq (0,0,0) \, {\rm and } \, \, \xi_1 x_1 + \xi_2 x_2 = z x_0 \, \, \\ & {\rm where} \, \, (x_0,x_1,x_2) \in \Omega \}
\end{align*}
which is also consonant with $\lambda_1 \sigma_2 - \lambda_1 \sigma_1 + \gamma$.
\smallskip

By \cite[Proposition 10.5.4]{KLMV}, $S$ is equal to the identity at infinity if and only if $W(\xi_1,\xi_2,z)$ 
is given by \eqref{eqR1}.
Then, by Theorem \ref{realizationPreLem}, Assumption \eqref{realizationThB} in Theorem \ref{realizationTh} is essentially equivalent to Assumption \eqref{JCF_classificationItem2} in Theorem \ref{JCF_classification}.
It remains to use Theorem \ref{JCF_classification} to complete the proof.
\end{pf}

%%%%%%%%%%%%%%%%%%%%%%%%%%%%%%%%%%%%%
%									%
%	NJCF							%
%									%
%%%%%%%%%%%%%%%%%%%%%%%%%%%%%%%%%%%%%
\subsection{Characterization theorem of the normalized joint characteristic functions}

Using the theory of algebraic curves there exists a compact Riemann surface $X$ together with a holomorphic mapping $\pi : X \rightarrow \mathbb P ^2$
such that $\pi$ is injective on the inverse image of the set of smooth points of $C$.
Then the kernel bundle $E$ (and similarly, $\widetilde{E}$) becomes a vector bundle over the inverse of $\pi$ along the non-singular points of $C$.
\smallskip

In order to continue, we add two assumptions on the discriminant curve $C$.
The first, is to assume that there exists an irreducible polynomial $f$
such that $p(\lambda) =f(\lambda)^r$. 
This makes $C$ an irreducible curve of degree $m=M/r$.
The second, is the assumption that the determinantal representation $U(\lambda)$ and $\widetilde{U}(\lambda)$ are maximal.
By definition, maximality of the determinantal representations means that for any $\lambda \in C$, 
$E(\lambda)$ is of maximal possible dimension, i.e. ${\rm dim} E(\lambda) = sr$, 
where $s$ is the multiplicity of the point $\lambda$ on $C$.
The maximality assumption ensures us, and this is a crucial point, that the vector bundle over $\pi^{-1}(C_{{\text smooth}})$
can be extended to a vector bundle of rank $r$ over $X$ (see \cite[Theorem 2.1]{MR97m:30051}).
\smallskip

The vector bundles associated to the input and output determinantal representations, 
denoted again by $E$ and $\widetilde{E}$, are flat vector bundles $V_\chi$ and $V_{\widetilde{\chi}}$,
where $\chi, \widetilde{\chi} :\pi_1(X) \rightarrow GL(r,\mathbb C)$ are factors of automorphy for the 
group of deck transformations on the universal covering $\widetilde{X}$ of $X$.
The vector bundle $V_\chi$ satisfies $h^0(V_\chi \otimes \Delta)=0$,
that is, the
twist $V_\chi \otimes \Delta$ of $V_\chi$ with a bundle $\Delta$ of half-order differentials on X (i.e.,
$\Delta \otimes \Delta = K$ is the line bundle of holomorphic differentials on $X$) has no global holomorphic sections.
\smallskip

It is shown in \cite{MR97m:30051, MR1704479} how to construct $M\times M$ matrices
($M = mr$) such that the map
\[
f(p) \rightarrow \frac{1}{\omega(p)} f(p) u^{\times}(p),
\]
(where $f(p)$ is a holomorphic half-ordered differential on $X$,
$\omega(p)$ is a meromorphic differential on $X$ and) 
establishes an isomorphism
from $V_\chi$ to $E\otimes\Delta\otimes O(1)$, where $E$ is the kernel bundle associated with
the maximal determinantal representation $U(\lambda) = \lambda_0 \gamma+\lambda_1 \sigma_2-\lambda_2 \sigma_1$ for
the curve $C$. Here $u^\times(p)$ is a {\it matrix of normalized sections} for $E$, 
i.e., a multiplicative $\Delta$-valued frame for $E$,
normalized to have poles exactly at the points of $C$ at infinity.
\smallskip

Thus the mapping $S$ is lifted to a bundle map $T$ operating between the
bundles $V_{\chi}$ and $V_{\widetilde{\chi}}$, which is related to $S$ via the matrices of normalized
sections implementing the equivalence between $V_{\chi}$ and $E$ and between $V_{\widetilde{\chi}}$
and $\widetilde{E}$ by
\begin{equation}
\label{eqST}
S(p) u ^\times(p) = \widetilde{u} ^\times(p)T(p).
\end{equation}
Under this isomorphism, we have the following scalar product
\begin{equation}
\label{eqCh5A}
\innerProductIndefiniteQuad{u^{\times}(p)}{u^{\times}(q)}{p,q}{E} =
K(\chi; p,\overline{q})
\quad {\rm whenever} \quad p \neq \overline{q},
\end{equation}
and
\begin{equation}
\label{eqCh5B}
\innerProductIndefiniteQuad{u^{\times}(p)}{u^{\times}(p)}{p,\overline{p}}{E}  =
I,
\end{equation}
where $K(\chi; p,q)$ is the Cauchy kernel associated to $\chi$.
Thus, the JCF metric properties as given in \eqref{jcfCond1} and \eqref{jcfCond2} are translated, using \eqref{eqCh5A} and \eqref{eqCh5B}, into
$T(p) \overline{T(\overline{p})} = 1$
and the kernel
\begin{equation}
\label{eqKerT}
K(\widetilde{\chi}; p,\overline{q})
-
T(p)
K(\chi; p,\overline{q})
    T(q)^*
,
\end{equation}
has $\kappa$ negative squares.
\smallskip

Therefore, the characterization theorem of the class of normalized joint characteristic functions of vessels over Pontryagin spaces can be stated as in the following theorem.
\begin{theorem}
\label{realizationTh1}
Let $T(p)$ be a vector bundle mapping on a real compact Riemann surface $X$ corresponding to $\chi$ and $\widetilde{\chi}$.
Then $T(p)$ is the normalized joint characteristic function of a vessel $\mathcal V$ over Pontryagin space with $\kappa$ negative index and with discriminant polynomial $p(\lambda_1,\lambda_2)$ with maximal input and output determinantal representations corresponding to $\chi$ and $\widetilde{\chi}$
if and only if
$T(p)$ is a non-zero holomorphic function in the neighborhood of $C$ at infinity, meromorphic on $X \setminus X_{\mathbb R}$,
$T(p)(\overline{T(\overline{p})}) = 1$ and the kernel \eqref{eqKerT} has $\kappa$ negative squares.
\end{theorem}

%%%%%%%%%%%%%%%%%%%%%%%%%%%%%%%%%%%%%
%									%
%	Applications					%
%									%
%%%%%%%%%%%%%%%%%%%%%%%%%%%%%%%%%%%%%
\section{Functional models associated to Pontryagin space vessels}
\label{secFuncModel}
\setcounter{equation}{0}

We hereby present the de Branges-Rovnyak theory for Pontryagin spaces of analytic sections on real compact Riemann surfaces using Liv\v{s}ic's vessel theory over indefinite inner product spaces.
Furthermore, we note that the (definite) de Branges theory presented in \cite{AVP1} was considered in the line-bundle case,
hence we use this opportunity to consider the generalization to the vector bundle case.
\smallskip

Recall that 
an irreducible curve $C$ in $\mathbb P^2$ with a maximal determinantal representation 
$U(\lambda) = \lambda_0 \gamma+\lambda_1 \sigma_2-\lambda_2 \sigma_1$ 
and an associated kernel bundle $E$ (with $E(\lambda) = \ker U(\lambda)$)
is equivalent to a compact Riemann surface $X$ and a
flat vector bundle $V_\chi$ of rank $m$ satisfying $h^0(V_\chi \otimes \Delta)=0$.
In this section we consider the converse direction.
We start with a real compact Riemann surface $X$ 
and assume that the vector bundles $V_\chi$ and $V_{\widetilde{\chi}}$ satisfy $h^0 (V_\chi \otimes \Delta)=0$ 
and $h^0 (V_{\widetilde{\chi}} \otimes \Delta)=0$.
We choose a pair of meromorphic
functions $y_1(p)$ and $y_2(p)$ which generates $\mathcal M (X)$, the field 
of meromorphic functions on X. 
Then, the mapping $p \xrightarrow{\pi} (y_1(p), y_2(p))$ is a birational embedding of $X$ to the
algebraic curve $C$ (we also use the notation $C_{y_1,y_2}$) defined by the projective closure of 
the image of $X$ under $\pi$ in $\mathbb P ^2$ (see \cite[Theorem 5.1]{MR1704479}). 
\smallskip

Let $T(p)$ be a multiplicative function on $X$ with multipliers corresponding to $\chi$ and $\widetilde{\chi}$.
We say that $T(p)$ is  a $(\chi,\widetilde{\chi})$-mapping with $\kappa$ negative squares 
if $T(p) \overline{T(p^{\tau})} = 1$ and the kernel
\begin{equation}
\label{zetaCont}
K(\widetilde{\chi}; p,\overline{q})
-
T(p)
K(\chi; p,\overline{q})
T(q)^*
\end{equation}
has $\kappa$ negative squares.
Then, the model space corresponding to $T$, denoted by $\mathcal P (T)$, is the reproducing kernel Pontryagin space 
with the reproducing kernel \eqref{zetaCont}.
\smallskip

For $y(u)$, a real meromorphic function defined on $X$  with simple poles,
the corresponding model operator, $M^{y}$ \cite[Equation 3-3]{MR1634421},
defined on sections of the vector bundle
$V_{\chi}\otimes \Delta$,
which are analytic in neighborhoods of the poles of $y(u)$.
It is given by
\[
M^{y}f(u)
=
y(u)f(u) +
\sum_{m=1}^{n}{c_m f(p^{(m)}) K(\chi; p^{(m)},u)}
,
\]
where $n$ is the degree of $y(u)$, the points $p_1,...,p_n$ are the $n$ distinct poles of $y(u)$ and
$c_1,...,c_n$ are the residues of $y(u)$ at these poles (up to a sign).
Moreover, the operator $M^y$ is bounded in $\mathcal P (T)$ and for any pair of meromorphic functions $y_1$ and $y_2$ the operators $M^{y_1}$ and $M^{y_2}$ commute.
\smallskip

Let us consider $R_{\alpha}^{y} = \left(M^{y} - \alpha I  \right)^{-1}$ 
where $\alpha$ is in the neighborhood of infinity.
Then we have (see \cite[Equation 3-4]{MR1634421})
\[
R_{\alpha}^{y}f(u) =
\frac{f(u)}{y(u) - \alpha} -
\sum_{j=1}^{n}{\frac{f(u^{(j)})}{dy(u^{(j)})}K(\chi; u^{(j)},u)},
\]
where $\left(u^{(j)}\right)_{j=1}^{n}$ is the set of $n$ points in $X$ such that $y(u^{(j)})=\alpha$.
It is the counterpart of the resolvent operator in the real compact Riemann surface case.
\smallskip

The (commutative) model vessel associated to a vessel $\mathcal V$ and a pair of real meromorphic functions is the collection
\begin{equation*}
\mathcal V_T =
\left(
 \, M^{y_1} \, , \,  M^{y_2}  \, ;  \, \mathcal P (T) \, , \,  \Phi \, ,  \, E  \, ;  \, \sigma_1 \, ,  \, \sigma_2 \,  , \,  \gamma  \, , \, \widetilde{\gamma}  \,
\right)
.
\end{equation*}
The vessels $\mathcal V_T$ and $\mathcal V$ are unitary equivalent on their principal subspaces and share the same normalized joint characteristic function.
The mapping between the inner space $\mathcal P$ of a vessel $\mathcal V_T$ to the model space is given by (see also \cite[Equation 3-5]{MR1634421}):
\begin{equation}
\label{ModelMap1}
p \rightarrow \frac{\xi_1 dy_1(z) + \xi_2 dy_2(z)}{\omega(z)}
P(\xi_1,\xi_2,z)\Phi(\xi_1A_1+\xi_2A_2 - \xi_1y_1(z) -\xi_2y_2(z))^{-1}p.
\end{equation}
Here $p \in \mathcal P$, $z \in X$, $\xi_1$ and $\xi_2$ are free parameters and $P$ is the projection of $E$ onto the output fiber $\widetilde{E} (p)$.

In \cite[Theorem 3.5]{AVP1}, we characterized the definite de Branges spaces defined on real compact Riemann surfaces.
The counterpart theory for the Pontryagin space setting is presented below.
\begin{theorem}
\label{preTh}
Let $X$ be a real compact Riemann surface.
Let ${\mathcal X}$ be a reproducing kernel Pontryagin space of sections of $V_{\widetilde{\chi}} \otimes \Delta$, 
with negative index $\kappa$, analytic in an open and connected set $\Omega$.
We choose real meromorphic functions $y_1$ and $y_2$ with simple poles generating $\mathcal M(X)$,
        such that $\Omega$ contains all the points above the singular points of $C_{y_1,y_2}$ and contains the poles of $y_1$ and $y_2$ and all the elements of $\mathcal X$ are regular at these points.
Furthermore, we assume that for every $\alpha,\beta \in \mathbb C $ such that their $n$ pre-images belong to $\Omega$, the following two conditions hold:
\begin{enumerate}[label=(\roman*)]
\item
${\mathcal X}$ is invariant under $R_\alpha^{y_1}$ and $R_\beta^{y_2}$.
\item
For every choice of $f,g \in {\mathcal X}$ such that $f$ and $g$ are analytic at the poles of $y_1$ and $y_2$,
the following identity
\begin{align}
\label{StructIdent4}
\innerProductIndefiniteReg{R_\alpha^{y_k}f}{g} & -
\innerProductIndefiniteReg{f}{R_\beta^{y_k}g} - (\alpha - \overline{\beta})
\innerProductIndefiniteReg{R_\alpha^{y_k}f}{R_\beta^{y_k}g}
= \nonumber
\\
&
-i(\alpha - \overline{\beta}) \sum_{l,t=1}^{n}
\frac{f(\nu^{(l)})}{dy_k(\nu^{(l)})}
K(\widetilde{\chi}; \nu^{(l)} , \overline{\omega^{(t)}})
\frac{\overline{g(\omega^{(t)})}}{\overline{dy_k(\omega^{(t)})}}
,
\end{align}
holds for $k=1,2$.
\end{enumerate}
        Then the reproducing kernel of $\mathcal X$ is of the form \eqref{zetaCont}
        where $T(\cdot)$ is a $(\chi,\widetilde{\chi})$-bundle mapping with $\kappa$ negative squares 
        for some flat unitary vector bundle $\chi$.
\end{theorem}

\begin{pf}
The complete proof in the Hilbert space setting is given in \cite[Theorem 3.4]{AVP1}.
The steps of the proof remain similar under the Pontryagin space and vector bundle generalizations.
\smallskip

We embed the commuting operators $M^{y_1}$ and $M^{y_2}$ in a commutative vessel over Pontryagin space
\[
( \, M^{y_1} \, ,  \,  M^{y_2} \, ; \, \mathcal X \, , \,  \Phi \, , \, \mathbb C^n \, ; \,  \sigma_1 \, , \,  \sigma_2 \, , \, \gamma  \, , \, \widetilde{\gamma} \, ),
\]
where
$\Phi: \mathcal{X} \rightarrow \mathbb C^n$ is the evaluation operator at infinity.
We define $\sigma_1$ and $\sigma_2$ by
\begin{equation}
\label{sigma12}
\sigma_k =
 \begin{psmallmatrix}
  c^1_k		& 				&  				&		&			&		$		$\\
  			&	\ddots 		&  				&  		&			&		$		$\\
  			&			 	& 	c^r_k		&  		&			&		$		$\\
  			&			 	& 				& 0		& c^{r+2}_k	&		&		&\\
  			&		 		&  				& c^{r+1}_k& 0		&		$		$\\
  			&		 		&  				&  		&			&\ddots	$		$\\
  			&		 		&  				&  		&			&		&0	& c^{r+2m}_k \\
  			&		 		&  				&  		&			&		& c^{r+2m-1}_k& 0\\
 \end{psmallmatrix}
.
\end{equation}
    where $c^1_k,...,c^r_k$ correspond to real poles of $y_k$,
    while $c^{r+1}_k,...,c^k_{r+2m}$ correspond to non-real poles 
    ($y_k$ is real and  hence non-real poles appear in conjugate pair).
Furthermore, let $\widetilde{\gamma}$ be defined by
\begin{equation*}
\widetilde{\gamma}_{j,k} =
\begin{cases}
c_{2}^{k} h_{1}^{k} - c_{1}^{k} h_{2}^{k}, 
& 
\overline{p^{(j)}} = p^{(k)} 
\\
\bracketsA{\overline{c_{2}^{j}} c_{1}^{k} - \overline{c_{1}^{j}} c_{2}^{k} }
    K(\widetilde{\chi} ; p^{(k)},\overline{p^{(j)}}),
& 
{\text otherwise},
\end{cases}
\end{equation*}
and $\gamma$ by
\begin{equation}
\label{input}
\gamma_{j,k} =
\Psi_{j,k}+
\begin{cases}
c_{2}^{k} h_{1}^{k} - c_{1}^{k} h_{2}^{k}
& 
\overline{p^{(j)}} = p^{(k)} 
\\
\bracketsA{\overline{c_{2}^{j}} c_{1}^{k} - \overline{c_{1}^{j}} c_{2}^{k} }
    K(\widetilde{\chi} ; p^{(k)},\overline{p^{(j)}})
,
& 
{\text otherwise},
\end{cases}
\end{equation}
where $\Psi_{j,k}= (\overline{c_1^j}c_2^k- \overline{c_2^j}c_1^k) K _{\mathcal{X}}(p^{(k)},\overline{p^{(j)}})$ 
and where $c_k^j$ and $h_k^j$ are the coefficients of $y_k$ at the pole $p^{(j)}$
\[
y_k(u) = - \frac{c^{m}_{k}}{t(u)}+h^{m}_{k} + o(|t|).
\]
Then, the colligation conditions \eqref{collCond} corresponding to $M^{y_1}$ and $M^{y_2}$ are equivalent to the structure identity for $M^{y_1}$ and $M^{y_2}$, respectively.
The proof of this claim is similar \cite[Lemma 5.7]{AVP1}, except that we use the collection formula as presented in \cite[Section 4]{av2}
instead the line-bundle version as presented in \cite[Lemma 4.1]{av3}.
\smallskip

Furthermore, using definitions (\ref{sigma12} -- \ref{input}), one can show that the output and the linkage (and hence the input) vessel conditions hold.
Also here, the proof of the generalizations to the indefinite and non--scalar case remains the same.
However, we note that $\gamma$ and $\widetilde{\gamma}$ are given in term of the non-scalar Cauchy kernel instead of the theta functions and the prime form.
\smallskip

The output determinantal representation is maximal since the canonical determinantal representation constructed from
a vector bundle and a pair of meromorphic functions is always maximal (see \cite[Theorem 5.1]{MR1704479}).
To show that the input determinantal representation is maximal uses same method as in \cite{AVP1}.
The point is that it is enough to show that all elements in $\Omega$ lie outside the joint spectrum of $M^{y_1}$ and $M^{y_2}$. 
Then, since the output determinantal representation is maximal, so does the input determinantal representation.
To prove this latter part, we associate to each $p_0 \in \Omega$ a new model vessel 
corresponding to a different pair of meromorphic functions $w_1$ and $w_2$ such that 
$(y_1(p_0),y_2(p_0))$ lies outside the joint spectrum of $M^{w_1}$ and $M^{w_2}$ and we conclude that $(y_1(p_0),y_2(p_0))$ 
lies also the joint spectrum of $M^{y_1}$ and $M^{y_2}$.
See Step 5 in the proof of \cite[Theorem 3.5]{AVP1} for the comprehensive details.
\smallskip

The mapping \eqref{ModelMap1} from $\mathcal P$ to the model space is the identity
and the proof in the vector bundle case is similar to step 5 in \cite{AVP1} with the following modifications.
The Cauchy kernels in the line bundle case are replaced by the non-scalar Cauchy kernels, using the collection formula as in \cite[Section 4]{av2}
and the normalized sections are matrices instead of vectors.
On the other hand, $\mathcal X$ is a reproducing kernel Pontryagin space and hence satisfies
\begin{equation}
\label{qwe123A}
h(z)
=
\innerProductIndefiniteReg{h}{K_{\mathcal X}(\cdot,z)}_{\mathcal X}.
\end{equation}
Therefore
\begin{align}
\nonumber
\hat{h}(z)
= &
\frac{\xi_1 dy_1(z) + \xi_2 dy_2(z)}{\omega(z)}
\widetilde{P}(\xi_1,\xi_2,z)\Phi \times
\\ &
\nonumber
\times
(\xi_1A_1+\xi_2A_2 - \xi_1 y_1(z) -\xi_2y_2(z))^{-1}h
\\ = &
\label{qwe123B}
h(z).
\end{align}
Combining \eqref{qwe123A} and \eqref{qwe123B}, we conclude that the reproducing kernel 
can be expressed explicitly in terms of the model space mapping by:
\begin{align*}
K_{\mathcal X}(p,q) =
&
\left(
\frac{\xi d y (p) }{w(p)}
\right)
\widetilde{P}(\xi,p)\Phi(\xi A - \xi y(p))^{-1}\times
\\
&
\nonumber
\hspace{2mm}\times
(\xi A - \xi y(q))^{-*}\Phi^*\widetilde{P}(\xi,q)^*
\overline{
\left(
\frac{\xi d y (q) }{w(q)}
\right)
}
.
\end{align*}
We use the relation (see Section $3$ in the preceding pages)
\begin{align*}
\nonumber
(\xi\sigma)
\Phi
(\xi A - \xi y(p))^{-1}
&
(\xi A - \xi y(q))^{-*}
\Phi^*
(\xi \sigma)^*
\\ & =
\frac
{ S(\xi,p) ( \xi \sigma ) S(\xi, q)^*  - ( \xi \sigma )}
{-i(\xi y(p) - \xi \overline{y(q)})}
,
\end{align*}
where $S(\xi,z)$ is given by
\[
S(\xi,z)
=
I - i ( \xi \sigma ) \Phi(\xi A - \xi y(z))^{-1}\Phi^*
\]
and we use \eqref{eqST} to conclude (see also \cite[Section 5]{AVP1}) the following
\begin{align*}
K_{\mathcal{X}}(p,q)
= &
K(\widetilde{\chi} ; p,q) -
T(p)
K(\chi;p,q)
T(q)^*.
\end{align*}
\end{pf}
%%%%%%%%%%%%%%%%%%%%%%%%%%%%%%%%%%%%%%%%%%%%%%%%%%%%%%%%%%%%%%%%%%%%%%%%%%%%%%%%%%%%%%%%%%%%%%%%%%%%%%%%%%%%%%%%%%%%%%%%

\begin{theorem}
\label{preThConv}
Let $X$ be a real compact Riemann surface,
let $T$ be a vector bundle mapping corresponds to $(\chi,\widetilde{\chi})$ with $\kappa$ negative squares
and let $\mathcal X$ be the associated reproducing kernel Pontryagin space with reproducing kernel of the form \eqref{zetaCont}.
Furthermore, let $y$ be a real meromorphic function on $X$ such that $T$ is regular at the poles of $y$.
Then, for any $\alpha \in \mathbb C$ such that all its pre-images under $y(\cdot)$ belong to $\Omega$, 
$\mathcal X$ is $R_\alpha^{y}$-invariant and the structure identity \eqref{StructIdent4} holds.
\end{theorem}

\begin{pf}
We choose an additional real meromorphic function $y_2$
such that all the poles of $y_2$ are contained in $\Omega$ and $y_1$ and $y_2$ generate $\mathcal M (X)$.
\smallskip

We associate to $\widetilde{\chi}$ and $\chi$ the canonical determinantal representations to ensure maximality of the input and output determinantal representations.
Then we apply the Pontryagin space realization theorem, i.e. Theorem \ref{realizationTh1}.
Hence there exists a two-operator commutative vessel over Pontraygin space
\[
\mathcal V = ( \, A_1 \, , \, A_2  \, ; \,  P  \, , \,  \Phi \, , \, \mathbb C ^n \, ; \,  \sigma_1 \, , \, \sigma_2 \, , \, \gamma \, , \,  \widetilde{\gamma} \, )
\]
such that $T$ is the normalized joint characteristic function of $\mathcal V$.
Here $P$ is an arbitrary Pontryagin space with negative index $\kappa$, $A_1$ and $A_2$ are bounded operators on $P$.
Then the associated model vessel is given by
\[
\mathcal V_T = ( \, M^{y_1} \, , \, M^{y_2}  \, ; \,  \mathcal P(T) \, , \, \Phi_T \, , \, \mathbb C ^n \, ; \,  \sigma_1 \, , \, \sigma_2 \, , \, \gamma \, , \,  \widetilde{\gamma} \, )
\]
and is an irreducible commutative two-operator vessel which is unitary equivalent,
on its principal subspace, to $\mathcal V$ and they share the same normalized joint characteristic function.
\smallskip

The space $\mathcal P(T)$ is invariant under $M^{y_1}$ and then, by \cite{av3}, is also $R_\alpha^{y_1}$-invariant for $\alpha$
in the neighborhood of infinity.
\smallskip

Finally, since the structure identity is equivalent to the colligation condition 
(Step 1 in the proof of \cite[Theorem 3.4]{AVP1}, the proof remains similar in this paper setting),
the colligation condition of $M^{y_1}$ in $\mathcal V_T$ implies the structure identity for $y_1(\cdot)$.
\end{pf}

%%%%%%%%%%%%%%%%%%%%%%%%%%%%%%%%%%%%
%									%
%	BL indefinte THM                %
%									%
%%%%%%%%%%%%%%%%%%%%%%%%%%%%%%%%%%%%
\section{Beurling theorem for indefinite Hardy spaces on Finite Bordered Riemann surfaces}

There were several attempts to study shift operators and generalize Beurling's Theorem to indefinite inner product Hardy spaces, 
see, in different contexts, \cite{MR543738} and the series of papers by Ball and Helton, of which we mention, \cite{MR695942} and \cite{MR756761}.
These papers consider the Hardy spaces of functions with values in $\mathbb C^n$ and use Halmos' wandering subspaces in order to 
state a Beurling's theorem on $\mathbb C ^n$ equipped with a Hermitian form.
\smallskip

One may replace the unit disk $\mathbb D$ by an arbitrary finite bordered Riemann surface
and study invariant subspaces of square summable sections of the vector bundle $V_\chi \otimes \Delta$ on a finite bordered Riemann surface
(as in the previous sections, $V_\chi$ is a flat vector bundle of rank $m$ and $\Delta$ is a square root of the canonical line bundle).
Let $S$ be a finite bordered Riemann surface with boundary $\partial S$ having $k>0$ components, denoted by $X_0,...,X_{k-1}$.
The double of $S$ is a compact Riemann surface $X$ with a natural antiholomorphic involution $\tau$
(the boundary $\partial S$ coincides with the set of fixed points of $\tau$ on $X$, denoted by $X_{\mathbb R}$)
turning $X$ into a real compact Riemann surface of dividing type.
Finally, we recall that a real meromorphic function of dividing type is a real meromorphic function on $X$ 
such that $p\in \partial S$ if and only if $y(p) \in \mathbb C_{\mathbb R}$ 
(the existence of such functions is proved in \cite{ahlfors}).
\smallskip

It is then natural to consider (see \cite{ad2})
the indefinite inner product given by associating an $m \times m$ signature matrix to each boundary component
(or a sign to each boundary component in the line bundle setting).
The indefinite inner product is given by
\[
    \innerProductIndefiniteTri{f}{g}{H^2_{\chi; J_0,...,J_{k-1}}} = 
\sum_{i=0}^{k-1}{ \int_{X_i}{g(p)^* J_i f(p)}},
\]
where $J_i$ is the signature matrix associated to $X_i$
(more generally, one may consider $J(\cdot)$ to be a continuous matrix function 
whose values are $m \times m$ selfadjoint nonsingular matrices)
and satisfying
\[
\chi(T)^* J_i(T \widetilde{p})\chi(T) = J_i(\widetilde{p})
\]
for any $\widetilde{p} \in  \widetilde{X_i}$, an element in the universal covering above the $i$-th boundary component,
and any deck transformation $T$ of $\widetilde{S}$ over $S$, which identified as element of $ \pi_1(S)$.
\smallskip

We assume the vector bundle $V_\chi \otimes \Delta$ satisfies $h^0(V_\chi \otimes \Delta )=0$.
Then the corresponding Hardy space, denoted by $H^2_{J_0,...,J_{k-1}}(S,V_\chi \otimes \Delta)$ 
(when no confusion is possible, we use the shortcut $H^2_{\chi; (J_{\ell})}$),
is nondegenerate see \cite[Proposition 2.2]{av2}.
Furthermore, $H^2_{\chi; (J_{\ell})}$ is a reproducing kernel Krein space.
\smallskip

We recall that in \cite[Section 4]{av2} the authors introduced an isometric isomorphism (up to a $2 \pi$ factor) 
from the indefinite Hardy space of sections defined on a finite bordered Riemann surface 
to the indefinite Hardy space on the unit disc $H^{2,M}_{\mathbf{J}}(\mathbb T)$ where $M = mn$
($m$ is the rank of the vector bundle $V_\chi$ while $n$ is degree of $y$)
and where the Hermitian form is given by
\[
    \innerProductIndefiniteTri{f}{g}{H^{2,M}_{\mathbf{J}}(\mathbb T)} = 
\frac{1}{2 \pi}\int_{0}^{2 \pi} g(e^{it})^* \mathbf{J} f(e^{it}) dt
,
\]
where $\mathbf{J}$ is an $M\times M$ selfadjoint matrix 
and $f$ and $g$ are $\mathbb C^M$-valued function on $\mathbb D$. 
The isomorphism is based on the boundary uniformization theorem of finite bordered Riemann surfaces and,
for some fixed $\lambda_0 \in \mathbb C$ with $n$ pre-images, the isometric isomorphism is given explicitly by 
\begin{equation}
    \label{eqIso}
    f(p)
    \rightarrow
    v_{i,f}(\lambda) = (\lambda - \lambda_0) \sum _{j=1}^{n}
    \frac{K(\chi ; \lambda_0^i,\lambda^j) f(\lambda^j) }
    {dz(\lambda^j) \, \sqrt{dz}(\lambda_0^i)}
    ,\quad
    i=1,...,n
\end{equation}
where $\lambda^i$ denotes the $i$-th pre-image of $\lambda \in \mathbb C$ under the uniformazation mapping $z(p)$.

\begin{lem}
\label{lemResolventInvReal}
Let $S$ be a finite bordered Riemann surface, 
let $y(\cdot)$ be a real meromorphic function of dividing type on the double $X$
and assume that $\mathcal H^\perp$ is a closed subspace in $H^2_{J_0,...,J_{k-1}}(S,V_\chi \otimes \Delta)$.
Then for any $\alpha \in \mathbb C _+$ the operator $\mathcal M _{\frac{1}{y(p)-\overline{\alpha}}}$ is bounded away from zero in norm.
\end{lem}

% Proof %
\begin{pf}
    We compute the image of the multiplication operator $\mathcal M _{\frac{1}{y(p)-\overline{\alpha}}} \, f$,
    where $f \in \mathcal H^\perp$ with $\normA{f}=1$, under the isomorphism \eqref{eqIso}.
    \smallskip

    We recall that the topology of $H^{2,M}_{\mathbf{J}}(\mathbb T)$ is, by definition, the 
    topology of $H^{2}_M(\mathbb T)$. Thus we have
\begin{align*}
    \normTreA{\frac{1}{2 \pi}\mathcal M _{\frac{1}{y(p)-\overline{\alpha}}}f(p)}{H^2_{\chi; (J_{\ell})}}{2}
    & =
    \normTreA{
        v_{\frac{f(\cdot)}{y(\cdot)-\overline{\alpha}}}(\lambda) }{H^{2,M}_{\mathbf{J}}(\mathbb T)}{2}
      =
    \normTreA{
        v_{\frac{f(\cdot)}{y(\cdot)-\overline{\alpha}}}(\lambda) }{H^{2,M}_{}(\mathbb T)}{2}
    \\ & =
    \int_{\mathbb T} 
    v_{\frac{f(\cdot)}{y(\cdot)-\overline{\alpha}}}(\lambda)^*  
    v_{\frac{f(\cdot)}{y(\cdot)-\overline{\alpha}}}(\lambda) 
    d \lambda 
    \\ & =
    \sum_{i=1}^{n}
    \int_{\mathbb T} 
    \overline{v_{i,\frac{f(\cdot)}{y(\cdot)-\overline{\alpha}}}(\lambda)}
    v_{i,\frac{f(\cdot)}{y(\cdot)-\overline{\alpha}}}(\lambda) 
    d \lambda. 
\end{align*}
By definition of the isomorphism, the last term becomes
\begin{align*}
    \sum_{i=0}^{n} \int_{\mathbb T} \abs{\lambda_0 - \lambda}^2 
    &
    \sum _{j=1}^{n}
    \overline{\frac{K(\chi ; \lambda_0^i,\lambda^j) }{dz(\lambda^j) \, \sqrt{dz}(\lambda_0^i)}}
    \frac{\overline{f(\lambda^j)}}{y(\lambda^j)- \alpha}
    \times \\ &
    \sum _{k=1}^{n}
    \frac{K(\chi ; \lambda_0^i,\lambda^k) }{dz(\lambda^k) \, \sqrt{dz}(\lambda_0^i)}    
    \frac{f(\lambda^k)}{y(\lambda^k)- \overline{\alpha}}
    \\ = 
    \int_{\mathbb T} \abs{\lambda_0 - \lambda}^2
    &
    \sum _{m,j=1}^{n}
    \frac{\overline{f(\lambda^j)}}{y(\lambda^j)- \alpha}
    \frac{f(\lambda^j)}{y(\lambda^j)- \overline{\alpha}}
    \times \\ &
    \left(
    \frac{1}{\overline{dz(\lambda^j)}}
    \sum_{i=0}^{n}  
    \frac{\overline{K(\chi ; \lambda_0^i,\lambda^j) } \, K(\chi ; \lambda_0^i,\lambda^k) }{dz(\lambda_0^i)}
    \frac{1}{dz(\lambda^k)} 
    \right).
\end{align*}
    Using the collection formula \cite[Section 3]{av2}
    \[
    \sum_{i=0}^{n}  
    \frac{K(\chi ;\overline{\lambda_1^j} \lambda_0^i) \, K(\chi ; \lambda_0^i,\lambda_2^k) }{dz(\lambda_0^i)}
    =
    \frac{\lambda_2 - \overline{\lambda_1}}{(\overline{\lambda_1} - \lambda_0)(\lambda_2 - \lambda_0)}
    K(\chi ;\overline{\lambda_1^j} \lambda_2^k)
    ,
    \]
    and recalling that $\lambda_0$ has $n$ distinct pre-images and $\lambda$ is real,
    the inner summation is zero whenever $m \neq j$.
    Furthermore, we recall that the inverse of the unit circle under the mapping $z$  is $\partial S$ 
    and since $y$ is a real dividing meromorphic function, 
    we have that $y(\lambda^j)\in \mathbb R$ along the integration path.
    Finally, combining the above with the assumption that $\imag{\alpha} > 0$, we conclude the following
\begin{align*}
    &
    \int_{\mathbb T} \abs{\lambda_0 - \lambda}^2
    \sum _{j=1}^{n}
    \frac{\abs{f(\lambda^j)}^2}{\abs{y(\lambda^j)-\real{\alpha}}^2 + \abs{\imag{ \alpha}}^2}
    \sum_{i=0}^{n}
    \overline{K(\lambda_0^i,\lambda^j)} 
    K(\lambda_0^i,\lambda^j)
    \\ \geq &
    \frac{1}{\abs{\imag{ \alpha}}^2}
    \int_{\mathbb T} \abs{\lambda_0 - \lambda}^2
    \sum _{j=1}^{n}
    \abs{f(\lambda^j)}^2
    \sum_{i=0}^{n}  
    \overline{K(\lambda_0^i,\lambda^j)} 
    K(\lambda_0^i,\lambda^j)
    \\ = & 
    \frac{1}{\abs{\imag{ \alpha}}^2}
    \normTreA{v_f(\lambda)}{H ^{2,M}_{\mathbf {J}}(\mathbb T)}{2}
    =  
    \frac{1}{\abs{\imag{ \alpha}}^2}
    \normTreA{\frac{f(p)}{2 \pi}}{H^2_{\chi; (J_{\ell})}}{2}
    > D_{\alpha}.
\end{align*}
Here $D_{\alpha}>0$ is independent of the choice of $f(p)$
and hence the multiplication operator $\mathcal M _{\frac{1}{y(u)-\overline{\alpha}}}$ is bounded away from zero.
\end{pf}

Therefore, repeating the same arguments as appeared in \cite[Lemma 4.10]{AVP1},
we may state the following result.

\begin{Cy}
    \label{MM1}
    Let $S$ be a finite bordered Riemann surface and $X$ is its double,
    let $y(\cdot)$ be a real meromorphic function of dividing type on $X$ 
    and let $\mathcal H^\perp$ be a closed subspace in $H^2_{\chi; (J_{\ell})}$.
    Furthermore, we assume that the elements of $\mathcal H^\perp$ have an analytic extension with bounded point evaluations 
    in a neighborhood of the poles of $y(\cdot)$
    Then the following hold:
    \begin{enumerate}
        \item 
            $R_{\alpha_0}^y$ is invariant and bounded in $\mathcal H^\perp$.
        \item
            The structure identity holds also for $\alpha_0$.
    \end{enumerate}
\end{Cy}
\begin{pf}
    It is a consequence of Lemma \ref{lemResolventInvReal} and repeating the proof in  \cite[Lemma 4.10]{AVP1}.
\end{pf}

We use de Branges structure theorem for Pontryagin spaces, i.e. Theorem \ref{preTh},
to give the following Beurling type theorem.
We highlight the additional assumption that $\mathcal H^\perp$ has a finite dimensional negative part
and hence Theorem \ref{preTh} can be applied.

\begin{Tm}
\label{thmIndefiniteBL}
    Let $S$ be a finite bordered Riemann surface such that $h^0(V_{\widetilde{\chi}} \otimes \Delta)=0 $,  
    let $X$ be its double and let
    $H^2_{\widetilde{\chi}; (J_{\ell})}$ be an associated indefinite Hardy space on $X_+$. 
    Let $y_1$ and $y_2$ be a pair of real meromorphic functions of dividing type generating $\mathcal M (X)$.
Furthermore, assume that the following conditions hold:
\begin{enumerate}
\item 
    $\mathcal H$ is a closed subspace of  $H^2_{\widetilde{\chi}; J_{\ell}}$ and is invariant under the multiplication operators 
    $\mathcal{M}_{\frac{1}{y_1(u)-\overline{\alpha}}}$ and $\mathcal{M}_{\frac{1}{y_2(u)-\overline{\alpha}}}$ for every $\alpha \in \mathbb C_+$.
\item 
    $\mathcal H^{\perp}$ is a non-degenerate subspace with finite dimensional negative part.
    Furthermore, the elements of $\mathcal H^{\perp}$ are assumed to have analytic extensions with bounded point evaluations 
    in neighborhoods of the poles of $y_1$ and $y_2$ and of the pre-images of the singular points of the curve $C_{y_1,y_2}$.
\end{enumerate}
Then there exists a $(\chi, \widetilde{\chi})$ bundle mapping $T$ with $\kappa$ negative squares such that
\[
    \mathcal H = T H^2_{\widetilde{\chi}; (J_\ell)}.
\]
\end{Tm}

The outline of this proof follows the proof in \cite[Theorem 4.3]{AVP1}.
Basically, the nonpositivity of the Hardy space does not influence significantly the proof outline, 
except a single delicate point, which has been discussed and elaborated above in Lemma \ref{lemResolventInvReal} and Corollary \ref{MM1}.
\smallskip

\begin{pf}[of Theorem \ref{thmIndefiniteBL}]
By assumption, $\mathcal H^{\perp}$,  the orthogonal complement subspace of $\mathcal H$,
is a non-degenerate subspace with finite dimensional negative part and hence  $\mathcal H^{\perp}$ is a Pontryagin space.
Furthermore, by assumption, $\mathcal H$ is invariant under the multiplication operators 
$\mathcal{M}_{\frac{1}{y_1(u)-\overline{\alpha}}}$ and $\mathcal{M}_{\frac{1}{y_2(u)-\overline{\beta}}}$ 
and thus $\mathcal H^{\perp}$ is invariant under the operators $R_{\alpha}^{y_1}$ and $R_{\beta}^{y_2}$ where $\alpha,\beta \in \mathbb C_+$.
This result was proved in \cite[Lemma 4.8]{AVP1} and remains valid in the indefinite setting.
\smallskip

Furthermore, we note that even in the indefinite setting, 
the Cauchy kernels are eigenfunctions of the resolvent operator and
the structure identity holds for resolvent operators corresponding to points in the upper half plane.
Hence, using Corollary \ref{MM1}, 
$\mathcal H ^\perp$ is invariant under the bounded operator $R_{\alpha}^{y_k}$ where $\alpha \in \mathbb R$ in the 
neighborhood of infinity and the structure identity can be extended to $\alpha \in \mathbb R$.
\smallskip

    Since $h^0(V_\chi \otimes \Delta)=0$ and using \cite[Proposition 2.2]{av2}, $X$ also satisfies $h^0(V_\chi \otimes \Delta)=0$ 
    and hence we consider the (unique) extension of vector bundles on $S$ to vector bundles on $X$ (see \cite[Proposition 2.1]{av2}).
    We then apply the indefinite version of the structure theorem, Theorem \ref{preTh}, 
    on $\mathcal H ^{\perp}$.
    Thus, $\mathcal H^{\perp}$ is a reproducing kernel Pontryagin space with the reproducing kernel
    \[
        K_{\mathcal H^{\perp}}(p,q) = {K(\widetilde{\chi};p,q) - T(p) K({\chi};p,q) T(q)^*},
    \]
where $T$ is a bundle mapping from $H^2_{\chi; (J_{\ell})}$ to 
    $H^2_{\widetilde{\chi}; (J_{\ell})}$ with $\kappa$ negative squares.
    Therefore, we can write
    $$\mathcal H ^{\perp} = H^2_{\widetilde{\chi}; (J_\ell)}  \ominus T H^2_{\chi; (J_\ell)}$$
    and so
    $\mathcal H =  TH^2_{\chi; (J_\ell)}$.
\end{pf}

%%%%%%%%%%%%%%%%%%%%%%%%%%%%%%%%%%%%%
%			Bibliography			%
%%%%%%%%%%%%%%%%%%%%%%%%%%%%%%%%%%%%%

\def\cfgrv#1{\ifmmode\setbox7\hbox{$\accent"5E#1$}\else
  \setbox7\hbox{\accent"5E#1}\penalty 10000\relax\fi\raise 1\ht7
  \hbox{\lower1.05ex\hbox to 1\wd7{\hss\accent"12\hss}}\penalty 10000
  \hskip-1\wd7\penalty 10000\box7} \def\cprime{$'$} \def\cprime{$'$}
  \def\cprime{$'$} \def\lfhook#1{\setbox0=\hbox{#1}{\ooalign{\hidewidth
  \lower1.5ex\hbox{'}\hidewidth\crcr\unhbox0}}} \def\cprime{$'$}
  \def\cprime{$'$} \def\cprime{$'$} \def\cprime{$'$} \def\cprime{$'$}
  \def\cprime{$'$}

\end{document}